\documentclass[preprint, 3p, 12pt]{elsarticle}
\usepackage{soul}
\usepackage{graphicx, graphics}
\usepackage{syntonly}
\usepackage{amssymb, amsmath, amsfonts}
\usepackage{subfigure}
\usepackage{multirow}
\usepackage{color}
\usepackage{bm}
\usepackage{natbib}

\newcommand{\br}{\mathbf r}
\newcommand{\bx}{\mathbf x}
\newcommand{\bu}{\mathbf u}

\newcommand{\bv}{\mathbf v}

\newcommand{\bF}{\mathbf F}
\newcommand{\bK}{\mathbf K}

\newcommand{\ud}{\text{d}}

\newcommand{\cS}{\mathcal S}
\newcommand{\tpartial}{\tilde{\partial}}
\newcommand{\tnabla}{\tilde{\nabla}}
\newcommand{\rank}[1]{\mbox{Rank}(#1)}
\newcommand{\sz}[1]{\mbox{Dim}(#1)}

\newcommand{\fracp}[2]{\frac{\partial {#1}}{\partial {#2}}}

\newcommand{\fF}{\mathcal F}

\newcommand{\bsig}{\boldsymbol \sigma}

\newcommand{\bp}{\mathbf p}

\begin{document}
\begin{frontmatter}
\title{Higher order nonlocal operator method}
\author[BU]{Huilong Ren}
\ead{Huilong.Ren@uni-weimar.de}
\author[HU,TU]{Xiaoying Zhuang\corref{cor1}}
\ead{zhuang@ikm.uni-hannover.de;}
\author[KU0,KU]{Timon Rabczuk\corref{cor2}}
\ead{timon.rabczuk@tdt.edu.vn}
\cortext[cor1]{Corresponding author}
\cortext[cor2]{Corresponding author}
\address[KU0]{Division of Computational Mechanics, Ton Duc Thang University, Ho Chi Minh City, Viet Nam}
\address[KU]{Faculty of Civil Engineering, Ton Duc Thang University, Ho Chi Minh City, Viet Nam}
\address[BU]{Institute of Structural Mechanics, Bauhaus-University Weimar, 99423 Weimar, Germany}
\address[HU]{Institute of Conitnuum Mechanics, Leibniz University Hannover, Hannover, Germany}
\address[TU]{State Key Laboratory of Disaster Reduction in Civil Engineering, College of Civil Engineering,Tongji University, Shanghai 200092, China}
\begin{abstract}
We extend the nonlocal operator method to higher order scheme by using a higher order Taylor series expansion of the unknown field. Such a higher order scheme improves the original nonlocal operator method proposed by the authors in [A nonlocal operator method for solving partial differential equations], which can only achieve one-order convergence. The higher order nonlocal operator method obtains all partial derivatives with specified maximal order simultaneously without resorting to shape functions. The functional based on the nonlocal operators converts the construction of residual and stiffness matrix into a series of matrix multiplication on the nonlocal operator matrix. Several numerical examples solved by strong form or weak form are presented to show the capabilities of this method.

\end{abstract}
\end{frontmatter}

\section{Introduction}\label{sec:introduction}

In the field of solving Partial Differential Equations (PDEs), methods can be generally divided into (semi-)analytical methods and numerical methods. Analytical methods include the method of separation of variables \cite{churchill1963fourier}, integral transforms \cite{bargmann1961hilbert}, Homotopy Analysis Method (HAM) \cite{liao2003beyond}, Variational Iteration Method (VIM) \cite{he1999variational} and so on. Analytical methods have advantages in finding the approximate/exact solutions but are often restricted to regular geometry domain. The numerical methods contain Rayleigh-Ritz method, Finite Difference Method(FDM), Finite Element Methods(FEMs), Meshless Methods(MMs), isogeometric analysis \cite{hughes2005isogeometric}, to just name a few. In finite element methods, the computation domain is meshed into discrete elements and the shape function defined on the element is used to interpolate the field value within the element. Meshless methods comprise many different formulations \cite{nguyen2008meshless,chen2017meshfree}, for example, Smoothed Particle Hydrodynamics (SPH) \cite{lucy1977numerical,gingold1977smoothed}, Element-Free Galerkin method (EFG)\cite{Belytschko1994}, Reproducing Kernel Particle Method (RKPM) \cite{LiuJunZhang1995} and so on. Finite element methods and most of the meshless methods interpolate the field value in the domain by means of shape functions, and the derivatives in PDEs are constructed from the derivatives of the shape functions. Different from the methods by interpolation technique, finite difference method expresses the partial derivatives with finite difference. However, the finite difference method is only applicable for domain with regular geometry. For a more complete review of the PDEs by numerical methods, we refer to \cite{tadmor2012review}.

When it comes to higher order PDEs in higher dimensional space, finite element method, meshless methods and finite difference method confront some problems. For finite element methods, the topology of element in higher dimensional space is complicated. Though the simplex element is valid in any dimensions, the representation of the topology and calculation of the shape functions and their partial derivatives are cumbersome. Other difficulties involve the numerical integration and the continuity required on the interface between adjoint elements. Nevertheless, some finite element schemes are developed for arbitrary order of derivative (i.e. \cite{droniou2017mixed,schedensack2016new,wu2017nonconforming}). For meshless methods based on the shape functions, there is no problem for the mesh construction in the higher dimensional space. However, the numerical integration in meshless methods requires a background mesh, which is the same as the finite element methods. What's worse, the calculation of higher order derivative of the shape function is very expensive. One method to circumvent the numerical integration in background mesh is the nodal integration, which however surfers the rank-deficiency problem. The finite difference method can construct higher order finite difference to replace the higher order partial derivatives, but the stencil becomes more complicated. Other problems with higher order PDEs in high dimensional space involve the complicated boundary conditions at different orders of derivatives, the proof on uniqueness, robust, stability of the solution.

The fundamental elements in PDEs are various partial differential operators of different orders. How to deal with these operators is the central topic of various numerical methods. FEMs and most meshless methods start from the shape function for interpolation, while the derivatives of shape function are used to represent the differential operators. Such process is expensive for higher order differential derivatives in higher dimensions. The difficulties to numerically describe  the differential operators arise from the locality of the operator, where the locality denotes the operator being defined at a point. To circumvent the difficulties arising from locality, Nonlocal Operator Method (NOM) was proposed by the authors \cite{ren2019nom}. NOM starts from the common differential operators such as gradient, curl, divergence and Hessian operators, to define the nonlocal gradient, nonlocal curl, nonlocal divergence and nonlocal Hessian operators by introducing the support with finite characteristic length. These nonlocal operators can be viewed as the generalization of the local operators. When the support degenerates to one point, the nonlocal operators recover the local operators. Unlike FEMs, meshless methods or finite difference method, NOM is a ``true'' meshless method and only requires the neighbor list in the support in order to construct the nonlocal derivatives. The low order nonlocal operators \cite{ren2019nom}  can solve low order (not more than 4th-order) PDEs, but not higher order PDEs.

The purpose of the paper is to develop a higher order nonlocal operator method for solving higher order PDEs of multiple fields in multiple spatial dimensions. The nonlocal operator method obtains a set of partial derivatives of different orders at once. Combining with weighed residual method and variational principles, nonlocal operator method establishes the residual and tangent stiffness matrix for PDEs by some matrix operation on common terms, operator matrix. In contrast with finite element method or meshless method with shape functions, the nonlocal operator method leads to the differential operators directly and adopts the nodal integration method. The remainder of the paper is outlined as follows. In section \ref{sec:nom}, the basic concepts such as support and dual-support, and the low order nonlocal operators are reviewed and then the higher order nonlocal operator method based on Taylor series expansion of multiple variables is developed. We define a special quadratic functional to derive the nonlocal strong form for a $2n$-order PDEs based on the nonlocal operators in section \ref{sec:sqf}. We give several numerical examples to demonstrate the capabilities of this method in solving PDEs by strong form in section \ref{sec:numexamples} and by weak form in section \ref{sec:nwf}. Finally, we conclude in section \ref{sec:con}.

\section{Nonlocal operator method}\label{sec:nom}
\subsection{Basic concepts}
\begin{figure}[htp]
 \centering
 \subfigure[]{
 \label{fig:Coord}
 \includegraphics[width=.3\textwidth]{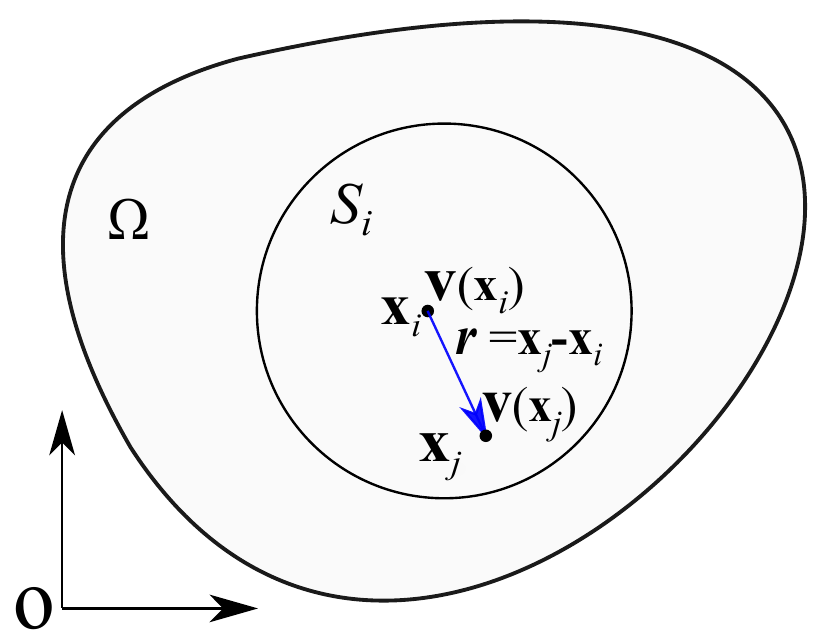}}
 \vspace{.1in}
 \subfigure[]{
 \label{fig:4support}
 \includegraphics[width=.35\textwidth]{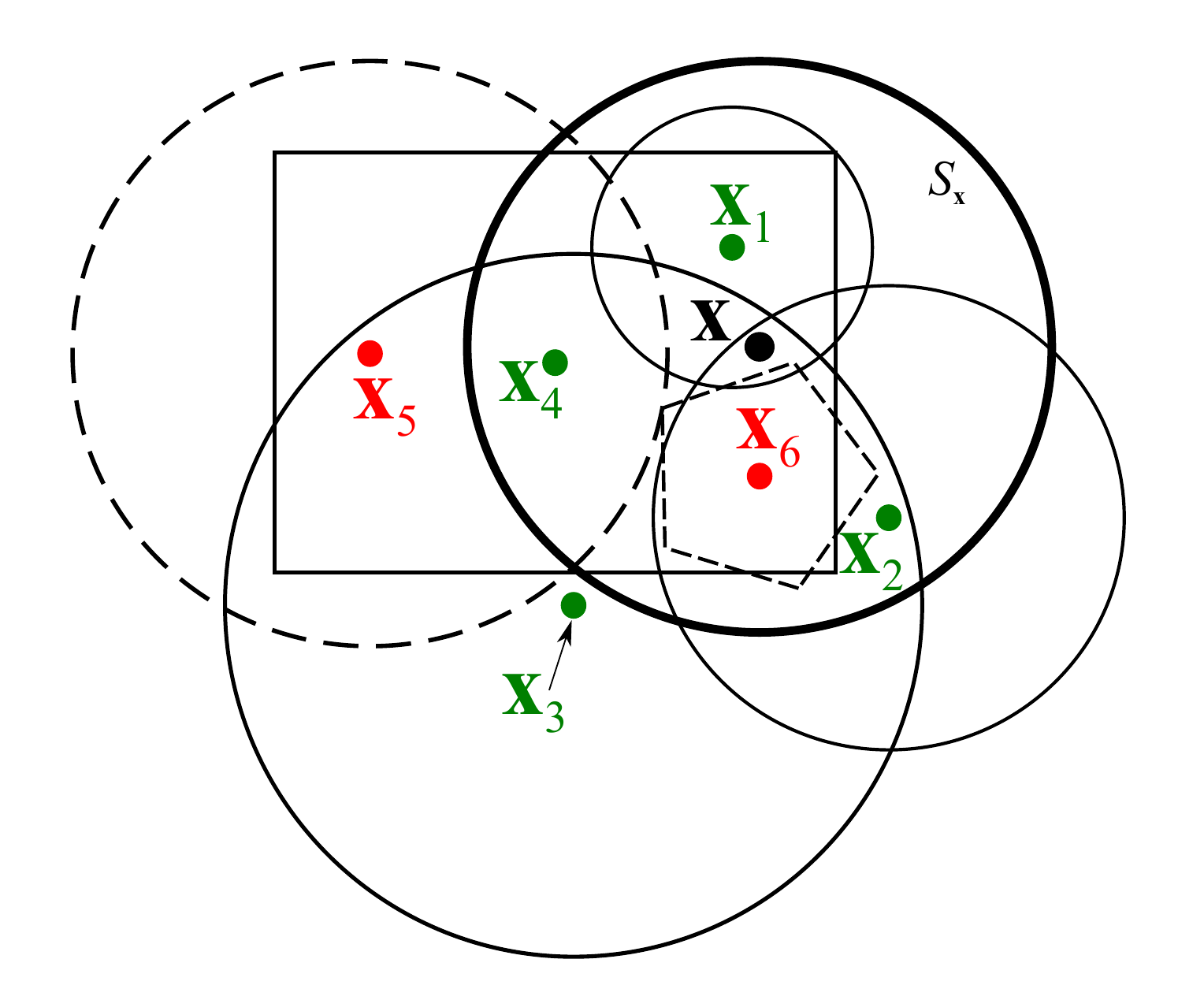}}\\
 \vspace{.3in}
\caption{(a) Domain and notation. (b) Schematic diagram for support and dual-support, all shapes above are supports, $\cS_{\bx}=\{\bx_1,\bx_2,\bx_4,\bx_6\} $, $\cS_{\bx}'=\{\bx_1,\bx_2,\bx_3,\bx_4\}$.}
\end{figure}
Consider a domain as shown in Fig.\ref{fig:Coord}, let $\bx_i$ be spatial coordinates in the domain $\mathbf \Omega$; $\br:=\bx_j-\bx_i$ is a spatial vector starts from $\bx_i$ to $\bx_j$; $\bv_i:=\bv(\bx_i,t)$ and $\bv_j:=\bv(\bx_j,t)$ are the field value for $\bx_i$ and $\bx_j$, respectively; $\bv_{ij}:=\bv_j-\bv_i$ is the relative field vector for spatial vector $\br$.

\textbf{Support} $\cS_{i}$ of point $\bx_i$ is the domain where any spatial point $\bx_j$ forms spatial vector $\br(=\bx_j-\bx_i)$ from $\bx_i$ to $\bx_j$. The support serves as the basis for the nonlocal operators. There is no restriction on the support shapes, which can be spherical domain, cube, semi-spherical domain, triangle and so on.

\textbf{Dual-support} is defined as a union of the points whose supports include $\bx$, denoted by 
\begin{align}
\cS_{i}'=\{\bx_j|\bx_i \in \cS_{j}\} \label{eq:dualsupport}.
\end{align}
Point $\bx_j$ forms dual-vector $\br'(=\bx_i-\bx_j=-\br)$ in $\cS_{i}'$. On the other hand, $\br'$ is the spatial vector formed in $\cS_{j}$. One example to illustrate the support and dual-support is shown in Fig.\ref{fig:4support}.

The nonlocal operator method uses the basic nonlocal operators to replace the local operator in calculus such as the gradient, divergence, curl and Hessian operators. The functional formulated by the local differential operator can be used to construct the residual or tangent stiffness matrix by replacing the local operator with the corresponding nonlocal operator. However, convergence rate of the original nonlocal operator is limited to 1 since the basic nonlocal operator is one-order. 

The nonlocal gradient of a vector field $\bv$ for point $\bx_i$ in support $\cS_{i}$ is defined as 
\begin{align}
\tnabla \bv_{i}:=\int_{\cS_{i}}w(\br) \bv_{ij} \otimes \br \ud V_{j} \cdot \Big(\int_{\cS_{i}}w(\br) \br \otimes \br \ud V_{j}\Big)^{-1} \label{eq:FGdef}.
\end{align}

The nonlocal gradient operator and its variation in discrete form are
\begin{align}
\tnabla \bv_{i}&=\sum_{j\in\cS_{i}}w(\br_j) \bv_{ij} \otimes \br_j \Delta V_{j} \cdot \Big(\sum_{j \in\cS_{i}}w(\br) \br \otimes \br \Delta V_{j}\Big)^{-1}  \label{eq:FGdef2},\\
\tnabla \delta\bv_{i}&=\sum_{j\in\cS_{i}}w(\br_j) \delta\bv_{ij} \otimes \br_j \Delta V_{j} \cdot \Big(\sum_{j\in\cS_{i}}w(\br) \br \otimes \br \Delta V_{j}\Big)^{-1}  \label{eq:FGdef3}.
\end{align}

The operator energy functional for vector field at point $\bx_i$ is
\begin{align}
\fF_i^{hg}&=p^{hg} \int_{\cS_i} w(\br) (\tnabla \bv_i\cdot\br-\bv_{ij})\cdot (\tnabla\bv_i\cdot\br-\bv_{ij}) \ud V_j
\end{align}
where $p^{hg}$ is the penalty coefficient. The residual and tangent stiffness matrix of $\fF_i^{hg}$ can be obtained with ease, we refer to \cite{ren2019nom} for more details.

\subsection{Higher order nonlocal operator method}\label{sec:honom}
Several formulations of the Taylor series expansion of a function of multiple variables are available in \ref{sec:Ty}. A scalar field $u_j$ at a point $j\in \cS_i$ can be obtained by the Taylor series expansion at $u_i$ in $d$ dimensions with maximal derivative order not more than $n$,
\begin{align}
u_j=u_i+\sum _{(n_{1},...,n_{d})\in \alpha_d^n}\frac {r_{1}^{n_{1}}...r_{d}^{n_{d}}}{n_{1}!...n_{d}!} u_{i,n_1...n_d}+O(r^{|\alpha|+1})\label{eq:mtse}
\end{align}
where
\begin{align}
\br&=(r_1,...,r_d)=(x_{j1}-x_{i1},...,x_{jd}-x_{id})\\
u_{i,n_1...n_d}&=\frac {\partial ^{n_{1}+...+n_{d}}u_i}{\partial x_{i1}^{n_{1}}...\partial x_{id}^{n_{d}}}\\
|\alpha|&=\max{(n_1+...+n_d)}
\end{align}
$\alpha_d^n$ is the list of flattened multi-indexes, where $d$ denotes the number of spatial dimensions and $n$ is the maximal order of partial derivative for one index. Two special multi-index can be written as
\begin{align}
\alpha_d^n&=\{(n_1,...,n_d)|1\leq\sum_{i=1}^d n_i\leq n,\, n_i\in \mathbb N^0, 1\leq i \leq d\}\label{eq:alpha2}\\
\mbox{or }\alpha_d^n&=\{(n_{1},...,n_{d})|1\leq \sum_{i=1}^d n_i,\, 0\leq n_i\leq n,n_i\in \mathbb N^0, 1\leq i\leq d\},\label{eq:alpha1}
\end{align}
where $\mathbb N^0=\{0,1,2,3,...\}$. Eq.\ref{eq:alpha1} gives a multi-index with $(1+n)^d-1$ elements and $|\alpha|=nd$, while the multi-index by Eq.\ref{eq:alpha2} has $\frac{(n+d)!}{n! d!}-1$ elements according to Combinatorics. In this paper, we adopt the multi-index by Eq.\ref{eq:alpha2} since it avoids the mixed higher order terms and has some benefit for numerical computation. The way to obtain all elements in $\alpha_d^n$ of Eq.\ref{eq:alpha2} by the Mathematica sees \ref{sec:mathcode}.

For any multi-index $(n_1,..., n_d)\in \alpha_d^n$, the partial derivative and the polynomial are
\begin{align}
u_{i,n_1...n_d},\, \frac {r_{1}^{n_{1}}...r_{d}^{n_{d}}}{n_{1}!...n_{d}!},\quad\forall (n_1,..., n_d)\in \alpha_d^n.
\end{align}

However, the original form of Taylor series expansion is very sensitive to the round-off error. For example, 
\[
{r_{1}^{n_{1}}...r_{d}^{n_{d}}}\propto h^{n_1+...+n_d}
\]
where $h$ is the characteristic length scale of the support. The higher order terms reduce to 0 quickly when $h<1$, or explode as $h>1$. It is expected to have the length scale $h$ approaching $1$.
When length scale of support $\cS_i$ at $u_i$ is taken into account, Taylor series expansion by Eq.\ref{eq:mtse} can be written as
\begin{align}
u_j&=u_i+\sum _{(n_1,...,n_d)\in \alpha_d^n}\frac {r_{1}^{n_{1}}... r_{d}^{n_{d}}}{h_i^{n_1+...+n_d}} \Big(\frac{h_i^{n_1+...+n_d}}{n_{1}!... n_{d}!}u_{i,n_1...n_d}\Big)+O(r^{n+1})\notag\\
&=u_i+\sum _{(n_1,...,n_d)\in \alpha_d^n}\frac {r_{1}^{n_{1}}... r_{d}^{n_{d}}}{h_i^{n_1+...+n_d}}\,\, u^h_{i,n_1...n_d}+O(r^{n+1})\label{eq:mtse2}
\end{align}
where $h_i$ is the characteristic length of $\cS_i$, and
\begin{align}
 u^h_{i,n_1...n_d}=\frac{h_i^{n_1+...+n_d}}{n_{1}!... n_{d}!}u_{i,n_1...n_d}
\end{align}

Let $\bp^h_j$, $\partial^h_\alpha u_i$ and $\partial_\alpha u_i$ be the list of the flattened polynomials, scaled partial derivatives, partial derivatives, respectively, based on multi-index notation $\alpha_d^n$ in Eq.\ref{eq:alpha2}, 
\begin{align}
\bp^h_j&=(\frac {r_{d}}{h},...,\frac {r_{1}^{n_{1}}... r_{d}^{n_{d}}}{h^{n_1+...+n_d}},...,\frac {r_1^n}{h^n})^T\label{eq:prh}\\
\partial^h_\alpha u_i&=(u^h_{i,0...1},...,u^h_{i,n_1...n_d},...,u^h_{i,n...0})^T\\
\partial_\alpha u_i&=(u_{i,0...1},...,u_{i,n_1...n_d},...,u_{i,n...0})^T.
\end{align}
Introducing $h$ in Eq.\ref{eq:prh} enables the terms in Eq.\ref{eq:prh} being in the ``same'' characteristic length scale.
The actual partial derivatives can be recovered by
\begin{align}
\partial_\alpha u_i=\mathbf H_i^{-1} \partial^h_\alpha u_i
\end{align}
where
\begin{align}
\mathbf H_i=\mbox{diag}\big[h_i,...,\frac{h_i^{n_1+...+n_d}}{n_1!... n_d!},...,\frac{h_i^{n}}{n!}\big]
\end{align}
where diag[$a_1, ..., a_n$] denotes a diagonal matrix whose diagonal entries starting in the upper left corner are $a_1, ..., a_n$.

Therefore, Taylor series expansion with $u_i$ being moved to left side of the equation can be written as
\begin{align}
u_{ij}=(\partial^h_\alpha u_i)^T\bp_j^h,\forall j\in \cS_i\label{eq:utse}
\end{align}
where $u_{ij}=u_j-u_i$.

Integrate $u_{ij}$ with weighted coefficient $w(\br) (\bp^h_j)^T$ in support $\cS_i$, we obtain
\begin{align}
\int_{\cS_i}w(\br) u_{ij} (\bp^h_j)^T \ud V_j&= (\partial^h_\alpha u_i)^T\,\int_{\cS_i} w(\br) \bp^h_j\otimes (\bp^h_j)^T \ud V_j \notag\\
&=(\partial_\alpha u_i)^T\,\mathbf H_i\,\int_{\cS_i} w(\br) \bp^h_j\otimes (\bp^h_j)^T \ud V_j 
\end{align}
where $w(\br)$ is the weight function. 

Therefore, the nonlocal operator $\tpartial_\alpha u_i$ can be obtained as
\begin{align}
\tpartial_\alpha u_i:=\mathbf H_i^{-1} \Big(\int_{\cS_i} w(\br) \bp^h_j\otimes (\bp^h_j)^T \ud V_j\Big)^{-1}\int_{\cS_i}w(\br) u_{ij} \bp^h_j \ud V_j=\bK_i \cdot\int_{\cS_i}w(\br) \bp^h_j u_{ij} \ud V_j \label{eq:tpui}
\end{align}
where 
\begin{align}
\bK_i:=\mathbf H_i^{-1} \Big(\int_{\cS_i} w(\br) \bp^h_j\otimes (\bp^h_j)^T \ud V_j\Big)^{-1}.\label{eq:Ki}
\end{align}

The reason to call Eq.\ref{eq:tpui} nonlocal operator is that it is defined in the support, in contrast with the local operator defined at a point. The nonlocal operator approximates the local operator with order up to $|\alpha|$. Traditional local operator is suitable for theoretical derivation but not for numerical analysis since its definition is limited to infinitesimal. The nonlocal operator can be viewed as a generalization of the conventional local operator.

The variation of $\tpartial_\alpha u_i$ is
\begin{align}
\tpartial_\alpha \delta u_i:=\bK_i\cdot\int_{\cS_i}w(\br) \bp^h_j (\delta u_{j}-\delta u_i) \ud V_j
\end{align}
In the continuous form, the number of dimensions of $\partial \delta u_i$ is infinite and discretization is required. 
After discretization of the domain by particles, the whole domain is represented by
\begin{align}
\Omega=\sum_{i=1}^{N} \Delta V_i
\end{align}
where $i$ is the global index of volume $\Delta V_i$, $N$ is the number of particles in $\Omega$.

Particles in $\cS_{i}$ are represented by 
\begin{align}
\cS_i&=\{j_1,...,j_k,...,j_{n_i}\}
\end{align}
where $j_1,...,j_k,...,j_{n_{i}}$ are the global indexes of neighbors of particle $i$, $n_i$ is the number of neighbors of $i$ in $\cS_i$. 

The discrete form of Eq.\ref{eq:tpui} and its variation are
\begin{align}
\tpartial_\alpha u_i&=\bK_i\cdot\sum_{j\in\cS_i}u_{ij}w(\br_j) \bp^h_j \Delta V_j=\bK_i \bp^h_{wi} \Delta \bu_i\label{eq:no}\\
\tpartial_\alpha \delta u_i&=\bK_i\cdot\sum_{j\in\cS_i}\delta u_{ij}w(\br_j) \bp^h_j \Delta V_j=\bK_i \bp^h_{wi} \delta\Delta \bu_i\label{eq:vno}
\end{align}
where
\begin{align}
\bK_i&=\mathbf H_i^{-1}\Big(\sum_{j\in\cS_i}w(\br)\bp_{j}^h\otimes(\bp_{j}^h)^T \Delta V_{j}\Big)^{-1},\label{eq:ibK}\\
\bp^h_{wi}&=\Big(w(\br_{j_1}) \bp^h_{j_1}\Delta V_{j_1},...,w(\br_{j_{n_i}}) \bp^h_{j_{n_i}}\Delta V_{j_{n_i}}\Big)\label{eq:bpw}\\
\Delta \bu_i&=(u_{i j_1},...,u_{i j_k},...,u_{i j_{n_i}})^T\label{eq:dbu}
\end{align}
When the weight function $w(\br)$ is selected as the reciprocal of the volume, Eq.\ref{eq:ibK} and Eq.\ref{eq:bpw} can be simplified further. The nonlocal operator provides all the partial derivatives with maximal order for single index up to $n$. The set of derivatives in PDEs of real application is a subset of the nonlocal operator. It should be noted that when the number of points in support is the same as the length of multi-index $\alpha_d^n$ and the coefficient matrix from Eq.\ref{eq:utse} for all points in support is well conditioned, the nonlocal operator can be obtained directly by the inverse of the coefficient matrix. In this case, the nonlocal operator serves as an efficient way to obtain the higher order finite difference scheme.

Each term in $\tpartial_\alpha u_i$ corresponds to the row of $\bK_i \bp^h_{wi}$ multiplying $\Delta\bu_i$. Eq.\ref{eq:no} can be used to replace the differential operators in PDEs to form the algebraic equations. This way is through strong form of the PDEs. The other ways to solve the linear (nonlinear) PDEs are through the weak formulations (weighted residual method) or the variational formulations (i.e. \cite{ren2019nom}). In these cases, the variation of $\partial_\alpha u_i$ in Eq.\ref{eq:vno} is required.

Eq.\ref{eq:no} can be written more concisely as
\begin{align}
\tpartial_\alpha u_i=\bK_i \bp^h_{wi} \Delta \bu_i=\mathbf B_{\alpha i}\mathbf u_i 
\end{align}
with $\mathbf B_{\alpha i}$ being the operator matrix for point $i$ based on multi-index $\alpha_d^n$  
\begin{align}
\mathbf B_{\alpha i}&=\begin{bmatrix}-(1,\cdots,1)_{n_p } \mathbf K_i \mathbf p_{wi}^h\\\mathbf K_i \mathbf p_{wi}^h\end{bmatrix}\\
\mathbf u_i&=(u_i,u_{j_1},u_{j_2},\cdots,u_{j_{n_i}})^T\end{align}
where $(1,\cdots,1)_{n_p} \mathbf K_i \mathbf p_{wi}^h$ is the column sum of $\mathbf K_i \mathbf p_{wi}^h$, $n_p$ is the length of $\alpha_d^n$. The operator matrix obtains all the partial derivatives of maximal order less than $|\alpha|+1$ by the nodal values in support. For real applications, one can select the specific rows in the operator matrix based on the partial derivatives contained in the specific PDEs. The template acts as

The traditional differential operator and their combination of one order or higher order and the corresponding variations can be constructed from Eq.\ref{eq:no} and Eq.\ref{eq:vno}, respectively. For example, the multi-index, polynomials and partial derivatives in two dimensions with maximal second-order derivatives are
\begin{align}
\alpha_2^2=&(01,02,10,11,20)\notag\\
\bp^h_j=&(y/h,y^2/h^2,x/h,x y/h^2,x^2/h^2)^T\notag\\
\tpartial_\alpha u_i=&(u_{,01},u_{,02},u_{,10},u_{,11},u_{,20})^T \label{eq:pu2d}
\end{align}
For the case of Poisson equation in 2D, $\nabla^2 u=f$. In the strong form, the operator $\nabla^2 u=\frac{\partial^2 u}{\partial x^2}+\frac{\partial^2 u}{\partial y^2}$ is required, one can select the $\partial_\alpha u_i[2]$ in Eq.\ref{eq:pu2d} for $\frac{\partial^2 u}{\partial y^2}$ and $\partial_\alpha u_i[5]$ in Eq.\ref{eq:pu2d} for $\frac{\partial^2 u}{\partial x^2}$. When solved in weak form, one can select the $\partial_\alpha u_i[1]$ in Eq.\ref{eq:pu2d} for $\frac{\partial u}{\partial y}$ and $\partial_\alpha u_i[3]$ in Eq.\ref{eq:pu2d} for $\frac{\partial u}{\partial x}$ to construct the tangent stiffness matrix.

In fact, the nonlocal operator $\partial_\alpha u_i$ in discrete form can be obtained by least squares. Consider the weighted square sum of the Taylor series expansion in $\cS_i$,
\begin{align}
\fF_i(\bu)=&\sum_{j \in\cS_i} w(\br)\big(u_{ij}-(\bp_j^h)^T\tpartial^h_\alpha u_i \big)^2 \Delta V_j\label{eq:Bu}\\
=&\sum_{j\in\cS_i} w(\br)\big(u_{ij}^2+\partial^h_\alpha u_i^T \bp^h_j(\bp_j^h)^T\partial^h_\alpha u_i-2 u_{ij} (\bp_j^h)^T\partial^h_\alpha u_i\big)\Delta V_j\notag\\
=&\sum_{j\in \cS_i} w(\br) u_{ij}^2\Delta V_j+\partial^h_\alpha u_i^T \sum_{j\in\cS_i}w(\br)\bp^h_j(\bp_j^h)^T \Delta V_j\,\,\partial^h_\alpha u_i-2 \Delta \bu_i^T\bp_{wi}^T \partial^h_\alpha u_i\notag\\
=&\sum_{j\in \cS_i} w(\br) u_{ij}^2\Delta V_j+\tpartial_\alpha u_i^T \mathbf H_i \sum_{j\in\cS_i}w(\br)\bp^h_j(\bp_j^h)^T \Delta V_j\,\,\mathbf H_i\partial^h_\alpha u_i-2 \Delta \bu_i^T(\bp_{wi}^{h})^T \mathbf H_i\tpartial_\alpha u_i\label{eq:fpu}
\end{align}
$\frac{\partial\fF_i(\bu)}{\partial (\tpartial_\alpha u_i)}=0$ leads to
\begin{align}
\tpartial_\alpha u_i=\mathbf H_i^{-1} \Big(\sum_{j\in\cS_i}w(\br)\bp^h_j(\bp_j^h)^T \Delta V_j\Big)^{-1}\bp_{wi}^h\Delta \bu_i =\bK_i \bp^h_{wi} \Delta \bu_i\label{eq:pu}
\end{align}
which is the same as Eq.\ref{eq:no}.

Meanwhile, Eq.\ref{eq:Bu} represents the operator energy functional in nonlocal operator method, and can be used to construct the tangent stiffness matrix of operator energy functional. The operator energy functional is the quadratic functional of the Taylor series expansion. Through Eq.\ref{eq:pu}, Eq.\ref{eq:fpu} can be simplified into
\begin{align}
\fF_i(\bu)
=&\sum_{j\in \cS_i} w(\br) u_{ij}^2\Delta V_j-\Delta \bu_i^T(\bp_{wi}^h)^T \Big(\sum_{j\in\cS_i}w(\br)\bp^h_j(\bp_j^h)^T \Delta V_j\Big)^{-1}\bp_{wi}^h\Delta \bu_i\notag\\
=&\Delta \bu_{i}^T \mathbf W_i \Delta \bu_i-\Delta \bu_i^T(\bp_{wi}^h)^T \Big(\sum_{j\in\cS_i}w(\br)\bp^h_j(\bp_j^h)^T \Delta V_j\Big)^{-1}\bp_{wi}^h \Delta \bu_i\notag\\
=&\Delta \bu_{i}^T \mathbf M_i \Delta \bu_i\label{eq:Bum}
\end{align}
where 
\begin{align}
\mathbf W_i&=\mbox{diag}\big[w(\br_{j_1})\Delta V_{j_1},...,w(\br_{j_{n_i}})\Delta V_{j_{n_i}}\big]\\
\mathbf M_i&=\mathbf W_i-(\bp_{wi}^h)^T \Big(\sum_{j\in\cS_i}w(\br)\bp^h_j(\bp_j^h)^T \Delta V_j\Big)^{-1}\bp_{wi}^h\label{eq:Mi}
\end{align}

The first and second variation of $\fF_i(\bu)$ read
\begin{align}
\delta\fF_i(\bu)&=2\Delta \bu_{i}^T \mathbf M_i \delta\Delta \bu_i.\label{eq:dBum}\\
\delta^2\fF_i(\bu)&=2\delta\Delta \bu_{i}^T \mathbf M_i\delta\Delta \bu_i.\label{eq:ddBum}
\end{align}

Let $\bv_i(j)=\sum_{k=1}^{n_i} \mathbf M_{i}(j,k)$ be the sum of row of matrix $\mathbf M_i$, then the tangent stiffness matrix can be extracted from Eq.\ref{eq:ddBum},
\begin{align}
\mathbf K_i^{hg}=\frac{p_{hg}}{m_i}\begin{bmatrix}\sum \bv_i& -\bv_i^T\\-\bv_i& \mathbf M_i \end{bmatrix}
\end{align}
where the first row(column) denotes the entries for point $i$, while the neighbors start from the second row(column), $p_{hg}$ is the penalty coefficient and $m_i$ the normalization coefficient
\begin{align}
m_i=\sum_{j\in\cS_i} w(\br) \br\cdot \br \Delta V_j,
\end{align}
where $\br$ varies for each $j$.

Let $n_i$ be the number of neighbors in $\cS_i$ and $n_p$ be the length of $\tpartial^h_\alpha u_i$. The dimensions of terms in $\mathbf M_i$ are
\begin{align}
&\sz{\mathbf W_i}=n_i\times n_i,\sz{\bp_{wi}^h}=n_p \times n_i,  \notag\\
&\sz{\Big(\sum_{j\in\cS_i}w(\br)\bp^h_j(\bp_j^h)^T \Delta V_j\Big)}=n_p\times n_p\notag\\ 
&\rank{(\bp_{wi}^h)^T \Big(\sum_{j\in\cS_i}w(\br)\bp^h_j(\bp_j^h)^T \Delta V_j\Big)^{-1}\bp_{wi}^h}\leq\min{(n_p,n_i)}.\notag 
\end{align}
When $n_i<n_p$, $\sum_{j\in\cS_i}w(\br)\bp^h_j(\bp_j^h)^T \Delta V_j$ is singular. It is required that $n_i\geq n_p$, so that ${\mathbf M_i}$ in Eq.\ref{eq:Mi} is well defined. The number of neighbors is selected as $5p+n_p$, where $p$ denotes the order of the nonlocal operator. These extra nodes are used to overcome the rank deficiency in the nodal integration.

The operator energy functional $\fF(\bu)$ represents the topology of the nonlocal operator method. Any field derived from $\fF(\bu)$ should try to satisfy $\fF(\bu)=0$ at the first step, which is independent with the actual physical model to be solved. 
\section{Quadratic functional}\label{sec:sqf}
A very special functional has the form 
\begin{align}
\fF=\frac{1}{2} \tpartial u^T \mathbf D \tpartial u\label{eq:sff}
\end{align}
where $\mathbf D$ is an arbitrary symmetric matrix, $\tpartial u \subset \tpartial_\alpha u$ in Eq.\ref{eq:no}. The operator matrix $\mathbf B$ is constructed from $\mathbf B_\alpha$ based on the index of terms $\partial u$ in $\partial_\alpha u$. Some examples of Eq.\ref{eq:sff} are given in \S \ref{sec:solidmechanics}.

When $\mathbf D$ is independent with the unknown functions $u$, the functional $\fF(\tpartial u)$ is pure quadratic, the first and second variation of $\fF(\tpartial u)$ at a point are
\begin{align}
\delta \fF&=\frac{\partial \fF}{\partial (\tpartial u)}=\tpartial \delta u^T \mathbf D \tpartial u=\delta \mathbf u^T \mathbf B^T \mathbf D\tpartial u\\
\delta^2 \fF&=\frac{\partial^2 \fF}{\partial (\tpartial u^T)\partial (\tpartial u)}=\tpartial \delta u^T \mathbf D \tpartial \delta u=\delta \mathbf u^T \mathbf B^T \mathbf D\mathbf B\delta \mathbf u
\end{align}
 and the residual and tangent stiffness matrix at a point can be written as
\begin{align}
\mathbf R(u)=\mathbf B^T \mathbf D\tpartial u, \,\mathbf K(u)=\mathbf B^T \mathbf D\mathbf B
\end{align}
When $\mathbf D:=\mathbf D (u,\partial u)$ is nonlinear tensor, the functional can be converted into quadratic functional by linearization and the Newton-Raphson can be employed to find the solution.

According to $\tpartial u_i\subset \tpartial_\alpha u_i$ in Eq.\ref{eq:tpui}, we write $\tpartial u_i=\mathbf K_i'\cdot \int_{\cS_i} w(\br) \mathbf p_j^h u_{ij} \ud V_j$, where $\mathbf K_i'\subset \mathbf K_i$ in Eq.\ref{eq:tpui}.   
Let $\bsig_i=\mathbf D\tpartial u_i$ and consider the variation in domain
\begin{align}
\delta\fF&=\delta\Big(\int_{\Omega} \frac{1}{2} \tpartial u_i^T \mathbf D \tpartial u_i \ud V_i\Big)
=\int_\Omega\,\tpartial u_i^T \mathbf D \tpartial \delta u_i\, \ud V_i\notag\\
&=\int_\Omega \bsig_i^T \tpartial \delta u_i \ud V_i
=\int_\Omega \bsig_i^T \bK_i' \int_{\cS_i} w(\br) \bp_j^h \delta u_{ij} \ud V_j \ud V_i\notag\\
&=\int_\Omega \int_{\cS_i} w(\br)\bsig_i^T \bK_i' \bp_j^h (\delta u_{j}-\delta u_{i}) \ud V_j \ud V_i\notag\\
&=\int_\Omega \Big(-\int_{\cS_i} w(\br)\bsig_i^T \bK_i' \bp_j^h \ud V_j+\int_{\cS'_i} \bsig_j^T \bK_j' \bp_i^h \ud V_j\Big)\delta u_{i} \ud V_i
\end{align}

Note that $\bp_i^h$ in $\cS'_i$ varies for different $j$ since $\bp_i^h$ is computed in $j$'s support $\cS_j$.

The terms with $\delta u_i$ in the first order variation $\delta \fF=0$ are
\begin{align}
-\int_{\cS_i} w(\br)\bsig_i^T \bK_i' \bp_j^h \ud V_j+\int_{\cS'_i} \bsig_j^T \bK_j' \bp_i^h \ud V_j,\label{eq:npde0}
\end{align}
with ``equivalent'' higher order partial differential term $-\partial^T (\mathbf D\partial u)$,
where $\partial:=(...,\frac {\partial ^{n_{1}+...+n_{d}}}{\partial x_{1}^{n_{1}}...\partial x_{d}^{n_{d}}},...)^T$ is the differential operator based on subset of multi-index $\alpha_d^n$ in Eq.\ref{eq:alpha2}. PDE given by $-\partial^T (\mathbf D\partial u)$ has a maximal differential order of $2n$. The nonlocal strong form by Eq.\ref{eq:npde0} can be solved directly by explicit integration algorithm. It should be noted that $\bsig, \tpartial u$ are in form of column vector for a scalar field $u$. The generalization of $u$ to vector field is straightforward.

Eq.\ref{eq:npde0} alone may suffer numerical instabilities (zero-energy mode), and therefore the operator energy functional by Eq.\ref{eq:Bu} is required. Eq.\ref{eq:npde0} with correction terms can be written as
\begin{align}
\tpartial_\alpha^T \bsig_i&\approx \int_{\cS_i}\big( w(\br)\bsig_i^T \bK_i' \bp_j^h+T^{hg}_{ij}\big) \ud V_j-\int_{\cS'_i} \big(\bsig_j^T \bK_j'\bp_i^h +T^{hg}_{ji}\big)d V_j\\
T^{hg}_{ij}&=w(\br)\frac{p^{hg}}{m_i}\big((\bp_j^h)^T\tpartial_\alpha^h u_i-u_{ij}\big)
\end{align}

\subsection{Elastic solid materials}\label{sec:solidmechanics}
In this section, we give some examples on how to express the linear/nonlinear elastic mechanics by the form of nonlocal operator method. The maximal derivative order in linear elastic mechanics is 2 and the corresponding weak form only requires first order partial derivative. 
The internal energy functional for plane stress, plane strain and 3D linear elastic solid at a point are 
\begin{align}
\fF_{plane\,stress}&=\frac{1}{2} \sigma:\varepsilon=\frac{1}{2}\tpartial \bu_{2d}^T D_{plane\,stress} \tpartial \bu_{2d}\\
\fF_{plane\,strain}&=\frac{1}{2} \sigma:\varepsilon=\frac{1}{2}\tpartial \bu_{2d}^T D_{plane\,strain} \tpartial \bu_{2d}\\
\fF_{3d}&=\frac{1}{2} \sigma:\varepsilon=\frac{1}{2}\tpartial \bu_{3d}^T D_{3d} \tpartial \bu_{3d},
\end{align}
where
\begin{align}
\tpartial \bu_{2d}&=(u_x,u_y,v_x,v_y)^T\\
\tpartial \bu_{3d}&=(u_x,u_y,u_z,v_x,v_y,v_z,w_x,w_y,w_z)^T
\end{align}
\begin{align}
D_{plane\, stress}=\frac{E}{1-\nu^2}\left[
\begin{array}{cccc}
 1 & 0 & 0 & \nu \\
 0 & \frac{1-\nu }{2} & \frac{1-\nu }{2} & 0 \\
 0 & \frac{1-\nu }{2} & \frac{1-\nu }{2} & 0 \\
 \nu & 0 & 0 & 1 \\
\end{array}
\right]
\end{align}
\begin{align}
D_{plane \, strain}=\frac{E}{(1-2\nu)(1+\nu)}\left[
\begin{array}{cccc}
 1-\nu & 0 & 0 & \nu \\
 0 & 1/2-\nu & 1/2-\nu & 0 \\
 0 & 1/2-\nu & 1/2-\nu & 0 \\
 \nu & 0 & 0 & 1-\nu \\
\end{array}
\right]
\end{align}
\begin{align}
D_{3D}=\left[
\begin{array}{ccccccccc}
 \lambda +2 \mu & 0 & 0 & 0 & \lambda & 0 & 0 & 0 & \lambda \\
 0 & \mu & 0 & \mu & 0 & 0 & 0 & 0 & 0 \\
 0 & 0 & \mu & 0 & 0 & 0 & \mu & 0 & 0 \\
 0 & \mu & 0 & \mu & 0 & 0 & 0 & 0 & 0 \\
 \lambda & 0 & 0 & 0 & \lambda +2 \mu & 0 & 0 & 0 & \lambda \\
 0 & 0 & 0 & 0 & 0 & \mu & 0 & \mu & 0 \\
 0 & 0 & \mu & 0 & 0 & 0 & \mu & 0 & 0 \\
 0 & 0 & 0 & 0 & 0 & \mu & 0 & \mu & 0 \\
 \lambda & 0 & 0 & 0 & \lambda & 0 & 0 & 0 & \lambda +2 \mu \\
\end{array}
\right]
\end{align}
The tangent stiffness matrix of that point can be extracted by performing the first or second order variation of the above functionals.

For nonlinear elastic material, the strain energy density is a function of the deformation gradient, i.e.
\[\fF(\bF)\]
while $\bF$ is consisted with the nonlocal operators in $\tpartial \bu_{3d}$
\begin{align}
\bF=\begin{bmatrix}F_1 & F_2 & F_3\\F_4& F_5 & F_6\\F_7& F_8 & F_9\end{bmatrix}=\begin{bmatrix} u_{x}+1 & u_{y} & u_{z} \\
 v_{x} & v_{y}+1 & v_{z} \\
 w_{x} & w_{y} & w_{z}+1 \\\end{bmatrix}.
\end{align}

Within the framework of total Lagrangian formulation, the first Piola-Kirchhoff stress is the direct derivative of the strain energy over the deformation gradient,
\begin{align}
\mathbf P=\frac{\partial \fF(\bF)}{\partial \bF}.
\end{align}
Furthermore, the material tensor (stress-strain relation) which is required in the implicit analysis can be obtained with the derivative of the first Piola-Kirchhoff stress,
\begin{align}
\mathbf D_4=\frac{\partial \mathbf P}{\partial \mathbf F}=\fracp{^2 \fF(\bF)}{\bF^T \partial \bF}.
\end{align}
The 4th order material tensor $\mathbf D_4$ can be expressed in matrix form when the deformation gradient is flattened.
\begin{align}
{\mathbf D}=\begin{bmatrix}
\frac {\partial P_1}{\partial F_{1}} & \frac {\partial P_1}{\partial F_{2}}&\cdots &\frac {\partial P_1}{\partial F_{9}}\\
\frac {\partial P_2}{\partial F_{1}} & \frac {\partial P_2}{\partial F_{2}}&\cdots &\frac {\partial P_2}{\partial F_{9}}\\
\vdots &\vdots &\ddots &\vdots \\
\frac {\partial P_9}{\partial F_{1}} & \frac {\partial P_1}{\partial F_{2}}&\cdots &\frac {\partial P_9}{\partial F_{9}}\\
\end{bmatrix}
=
\begin{bmatrix}
\frac {\partial ^{2}\fF(F)}{\partial F_{1}^{2}} & \frac {\partial ^{2}\fF(F)}{\partial F_{1}\,\partial F_{2}}&\cdots &\frac {\partial^{2}\fF(F)}{\partial F_{1}\,\partial F_{9}}\\
\frac {\partial ^{2}\fF(F)}{\partial F_{2}\partial F_{1}} & \frac {\partial ^{2}\fF(F)}{\partial F_{2}\,\partial F_{2}}&\cdots &\frac {\partial^{2}\fF(F)}{\partial F_{2}\,\partial F_{9}}\\
\vdots &\vdots &\ddots &\vdots \\
\frac {\partial ^{2}\fF(F)}{\partial F_{9}\partial F_{1}} & \frac {\partial ^{2}\fF(F)}{\partial F_{9}\,\partial F_{2}}&\cdots &\frac {\partial^{2}\fF(F)}{\partial F_{9}\,\partial F_{9}}\\
\end{bmatrix},
\end{align}
where the flattened deformation gradient and first Piola-Kirchhoff stress are
\begin{align}
F=(F_1, F_2, F_3,F_4,F_5, F_6,F_7,F_8,F_9)
\end{align}
and 
\begin{align}
P=\frac{\partial \fF(F)}{\partial F}=(\frac{\partial \fF(F)}{\partial F_1},\frac{\partial \fF(F)}{\partial F_2},\cdots,\frac{\partial \fF(F)}{\partial F_9}).
\end{align}

For the case of nearly incompressible Neo-Hooke material \cite{reese2000new}, the strain energy can be expressed as
\begin{align}
\fF(\bF)=\frac{1}{2} \kappa (J-1)^2+\frac{1}{2} \mu (\bF:\bF-3).\label{eq:nhm}
\end{align}
where $J=\det{\bF}$.

The first Piola-Kirchoff stress is
\begin{align}
\mathbf P=\frac{\partial \fF(\bF)}{\partial \bF}=\mu \bF+(J-1)\kappa J_{,\bF}.
\end{align}
With some derivation, the material tensor in matrix form can be written as
\begin{align}
\mathbf D=\mu \mathbf I_{9\times 9}+(J-1)\kappa J_{,F F}+\kappa J_{,F}\otimes J_{,F}
\end{align}
where $J_{,F}$ is the vector form of $J_{,\mathbf F}$, and 
\begin{align}J_{,F F}=
\left[
\begin{array}{ccccccccc}
 0 & 0 & 0 & 0 & F_{9} & \text{-}F_{8} & 0 & \text{-}F_{6} & F_{5} \\
 0 & 0 & 0 & \text{-}F_{9} & 0 & F_{7} & F_{6} & 0 & \text{-}F_{4} \\
 0 & 0 & 0 & F_{8} & \text{-}F_{7} & 0 & \text{-}F_{5} & F_{4} & 0 \\
 0 & \text{-}F_{9} & F_{8} & 0 & 0 & 0 & 0 & F_{3} & \text{-}F_{2} \\
 F_{9} & 0 & \text{-}F_{7} & 0 & 0 & 0 & \text{-}F_{3} & 0 & F_{1} \\
 \text{-}F_{8} & F_{7} & 0 & 0 & 0 & 0 & F_{2} & \text{-}F_{1} & 0 \\
 0 & F_{6} & \text{-}F_{5} & 0 & \text{-}F_{3} & F_{2} & 0 & 0 & 0 \\
 \text{-}F_{6} & 0 & F_{4} & F_{3} & 0 & \text{-}F_{1} & 0 & 0 & 0 \\
 F_{5} & \text{-}F_{4} & 0 & \text{-}F_{2} & F_{1} & 0 & 0 & 0 & 0 \\
\end{array}
\right].
\end{align}
The numerical example based on material model Eq.\ref{eq:nhm} is given in section \ref{sec:neoexample}.

\section{Numerical examples by strong form}\label{sec:numexamples}
The nonlocal operator defined in Eq.\ref{eq:no} can be used to replace the partial derivatives of different orders in the partial differential equation. In other word, we can use the nonlocal operator to solve the PDE by its strong form. In this sense, the nonlocal operator is similar to the finite difference method. However, finite difference scheme of different order is constructed on the regular grid, where the extension to higher dimensions or higher order derivative require special treatment, while the nonlocal operator is established simply based on the neighbor list in the support. In this section, we test the accuracy of nonlocal operator in solving second order ordinary differential equation (ODE) or PDE by strong form. Note that the operator energy functional is not required in solving PDE by strong form.

The first three numerical examples demonstrate the capabilities of nonlocal operator method in obtaining high order finite difference scheme. 
\subsection{Second-order ODE}
\begin{figure}
	\centering
		\includegraphics[width=9cm]{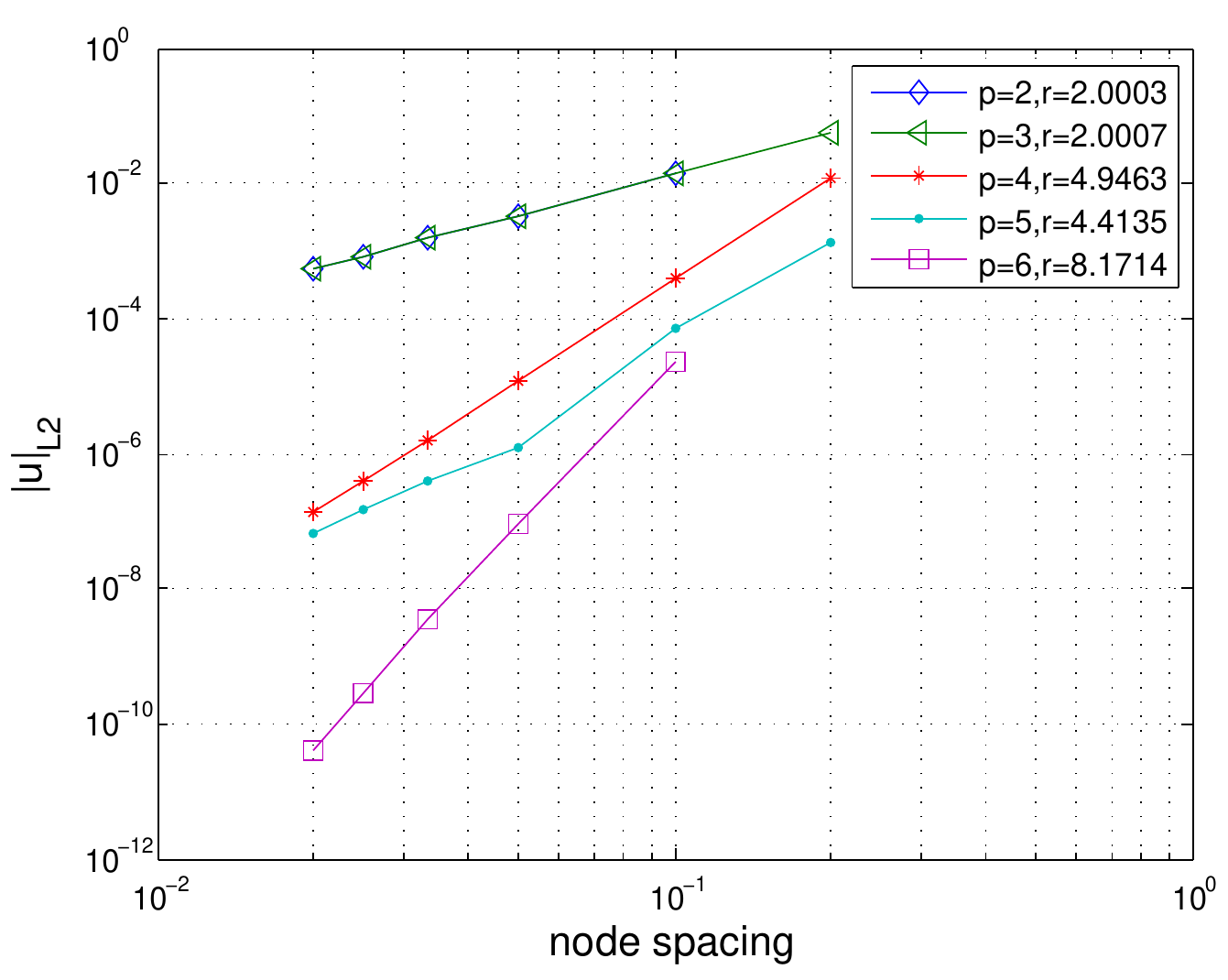}
	\caption{Convergence of the L2-norm for $u$.}
	\label{fig:barL2}
\end{figure}

The ODE with boundary condition is given by
\begin{align}
\frac{d ^2 u(x)}{d x^2}=20 x^3+\pi ^2 \cos (\pi x),u(0)=0,u(1)=0,x\in [0,1]
\end{align}
with analytic solution 
\[
u(x)=x^5-3 x-\cos (\pi x)+1.
\]
Since the highest order derivative in the ODE is two, the order of derivative in the nonlocal operator list should be $p\ge 2$. We test the nonlocal operator with $p=2,3,4,5,6$ in solving the second-order ODE. The minimal number of neighbors in the support is selected as the number of terms in the nonlocal operator. The difference between numerical result and theoretical solution is measured by 
the L2-norm, which is calculated by
\begin{align}
\|\bu\|_{L2}=\sqrt{\frac{\sum_j (\bu_j-\bu_j^{exact})\cdot (\bu_j-\bu_j^{exact}) \Delta V_j}{\sum_j \bu^{exact}_j\cdot \bu_j^{exact} \Delta V_j}}
\end{align}
The convergence of the L2-norm for $u$ is shown in Fig.\ref{fig:barL2}.
It can be seen that with the increase of order in the nonlocal operator, the convergence rate increases greatly. $p=2,3$ have the same convergence rate.
\subsection{1D Schr\"{o}dinger equation}
This section tests the accuracy of the eigenvalue problem in 1D. The Schr\"{o}dinger equation written in adimensional units for a one-dimensional harmonic oscillator is
\begin{align}
\big[-\frac{1}{2} \frac{\partial^2}{\partial x^2}+V(x)\big] \phi(x)=\lambda \phi(x),\quad V(x)=\frac{1}{2}\omega^2 x^2
\end{align}
For simplicity, we use $\omega=1$. The particles are uniformly distributed with constant spacing $\Delta x$ on the region [-10,10].

\begin{figure}
 \centering
 \subfigure[regular node distribution]{
 \label{fig:scheig1}
 \includegraphics[width=.45\textwidth]{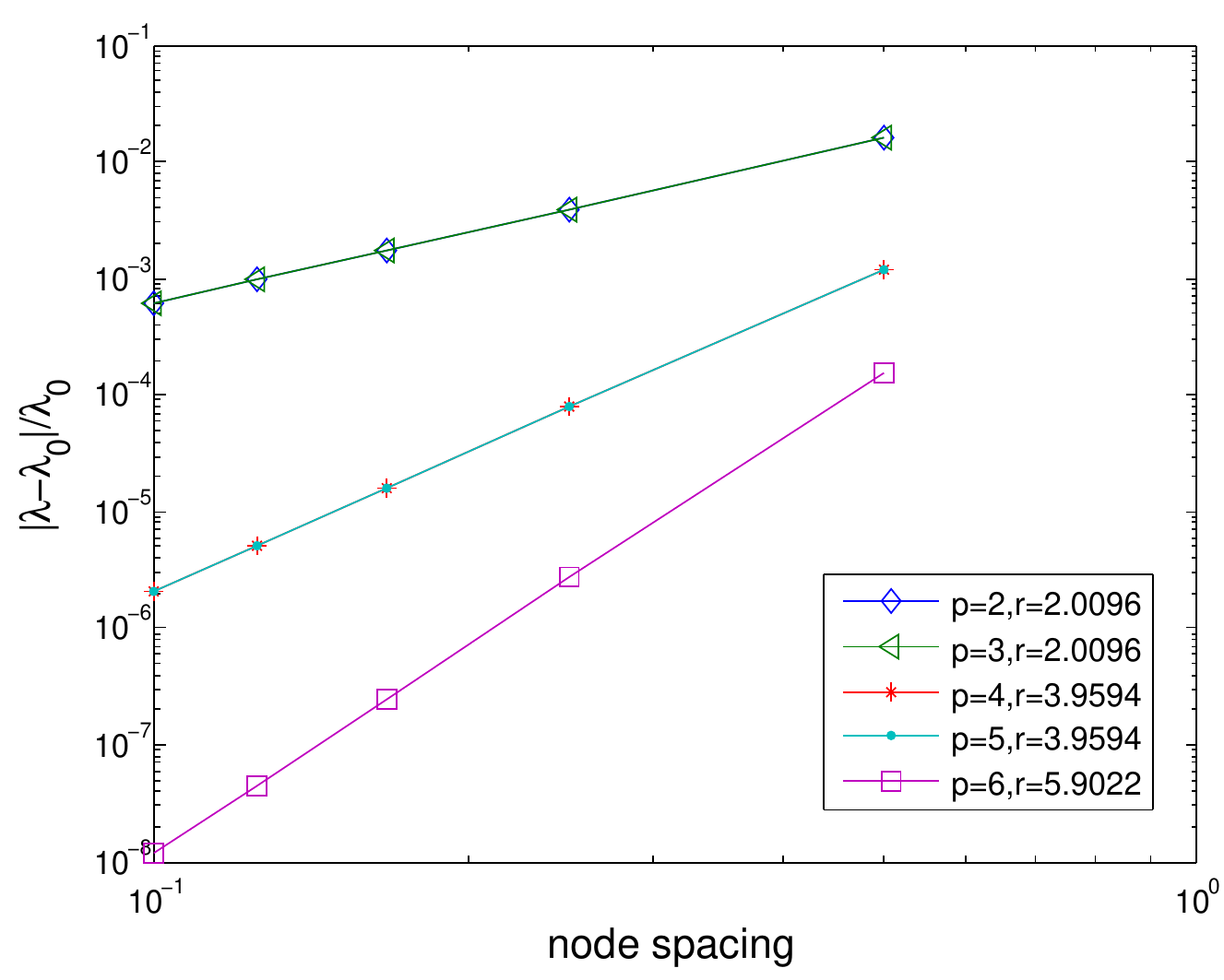}}
 \vspace{.1in}
 \subfigure[irregular node distribution]{
 \label{fig:scheig2}
 \includegraphics[width=.45\textwidth]{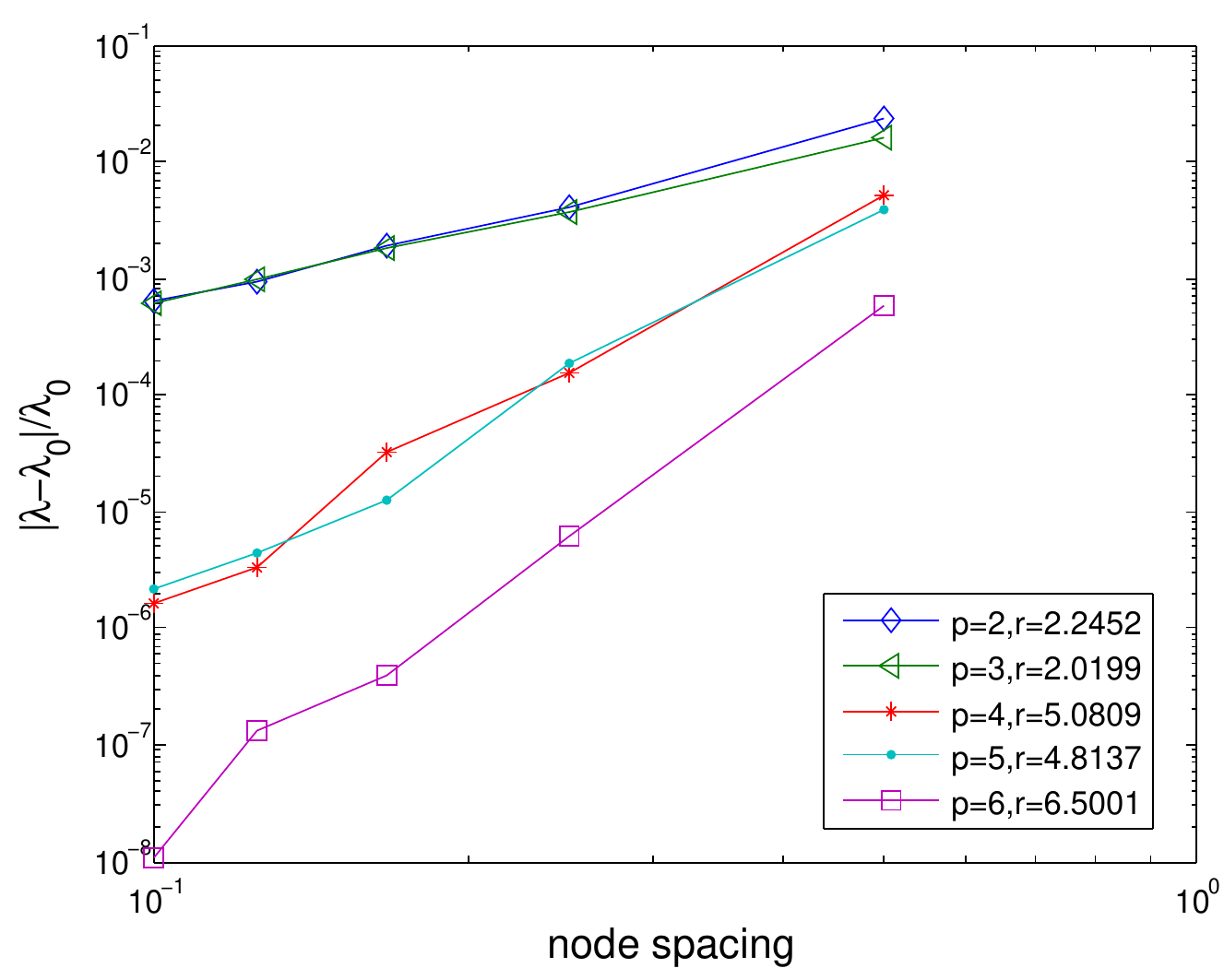}}
 \vspace{.3in}
\vspace{.1in}
\caption{Convergence of the lowest eigenvalue for a one-dimensional harmonic oscillator.}
\label{fig:scheig}
\end{figure}

The exact wave functions and eigenvalues can be expressed as
\begin{align}
\phi_n(x)=H_n(x) \exp(\pm \frac{x^2}{2}), \quad \lambda_n=n+\frac{1}{2}
\end{align}
where $n$ is a non-negative integer. $H_n(x)$ is the $n$-order Hermite polynomial.
We calculate the lowest eigenvalue and compare the numerical result with $\lambda_0=0.5$. The convergence plot of the error is shown in Fig.\ref{fig:scheig}.

\subsection{Poisson equation}
In this section, we test the Poisson equation 
\begin{align}
\nabla^2 u=2 x (y-1)(y-2 x+x y+2) e^{x-y},\quad (x,y)\in (0,1)\times (0,1)
\end{align}
with the boundary conditions
\begin{align}
u(x,0)=u(x,1)=0,\, x\in [0,1]\notag\\
u(0,y)=u(1,y)=0,\, y\in [0,1].\notag
\end{align}
The analytic solution is 
\begin{align}
u(x,y)=x(1-x)y(1-y)e^{x-y}.
\end{align}

\begin{figure}
 \centering
 \subfigure[]{
 \label{fig:mesh5}
 \includegraphics[width=.3\textwidth]{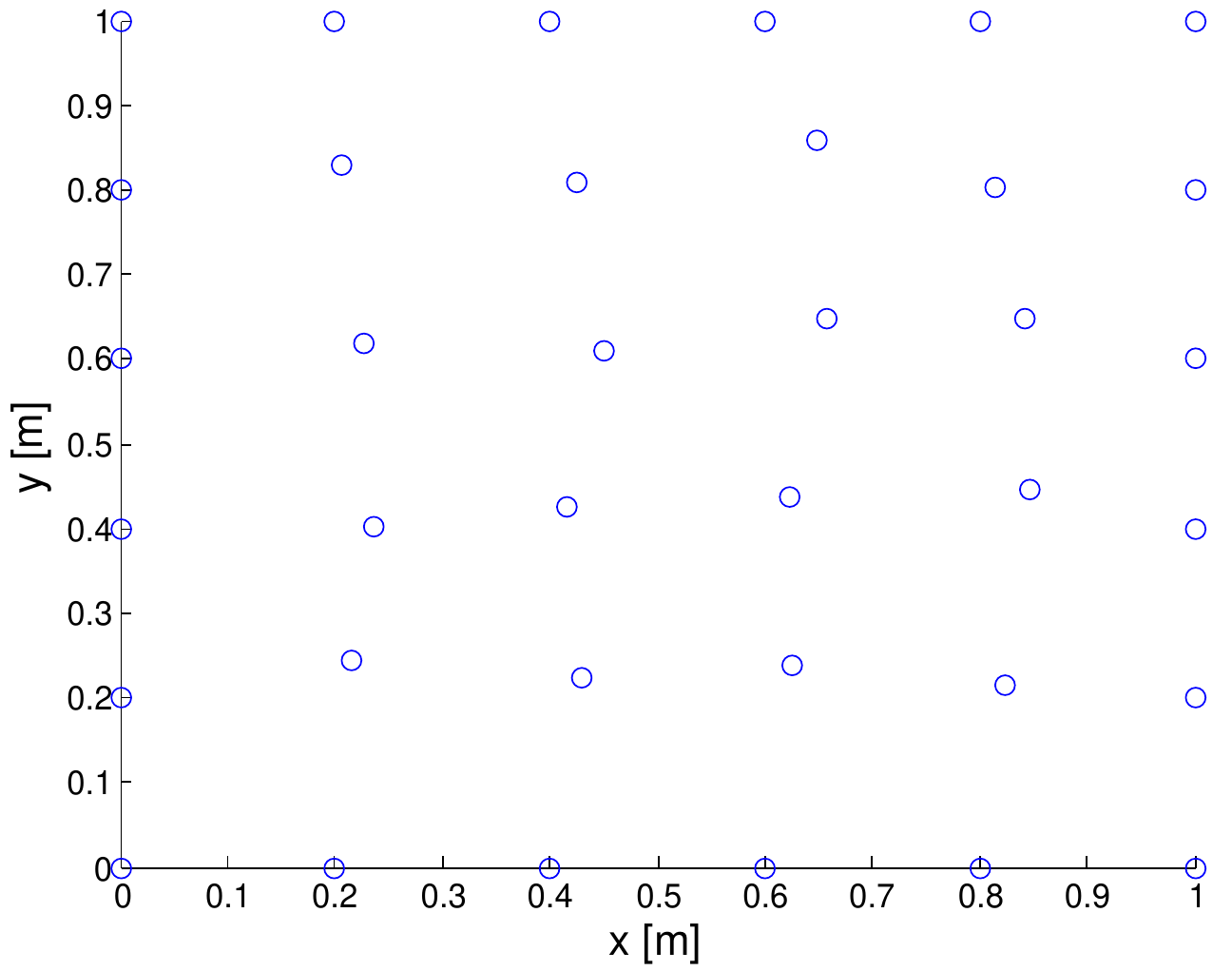}}
 \subfigure[]{
 \label{fig:mesh10}
 \includegraphics[width=.3\textwidth]{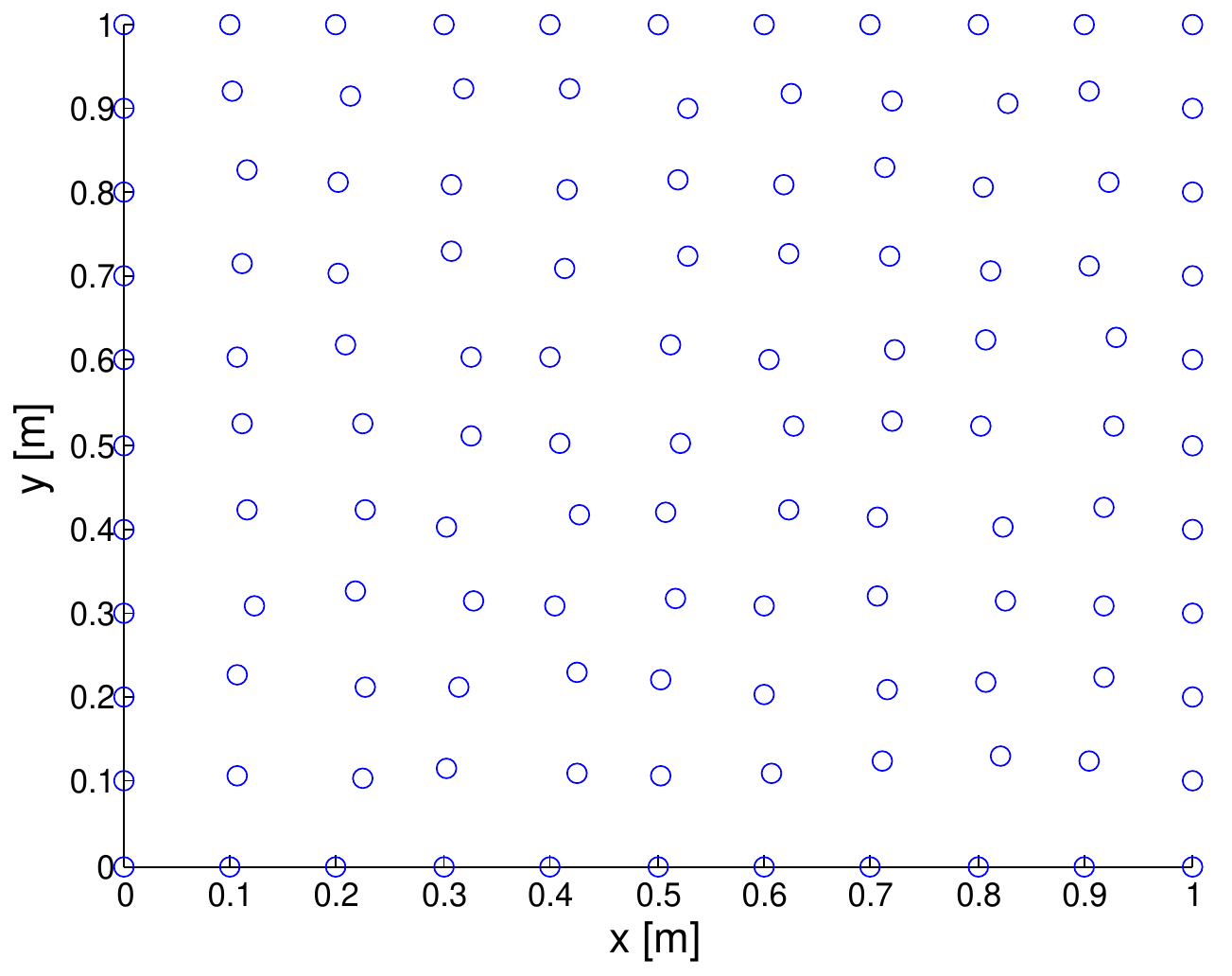}}
\vspace{.1in}
 \subfigure[]{
 \label{fig:mesh20}
 \includegraphics[width=.3\textwidth]{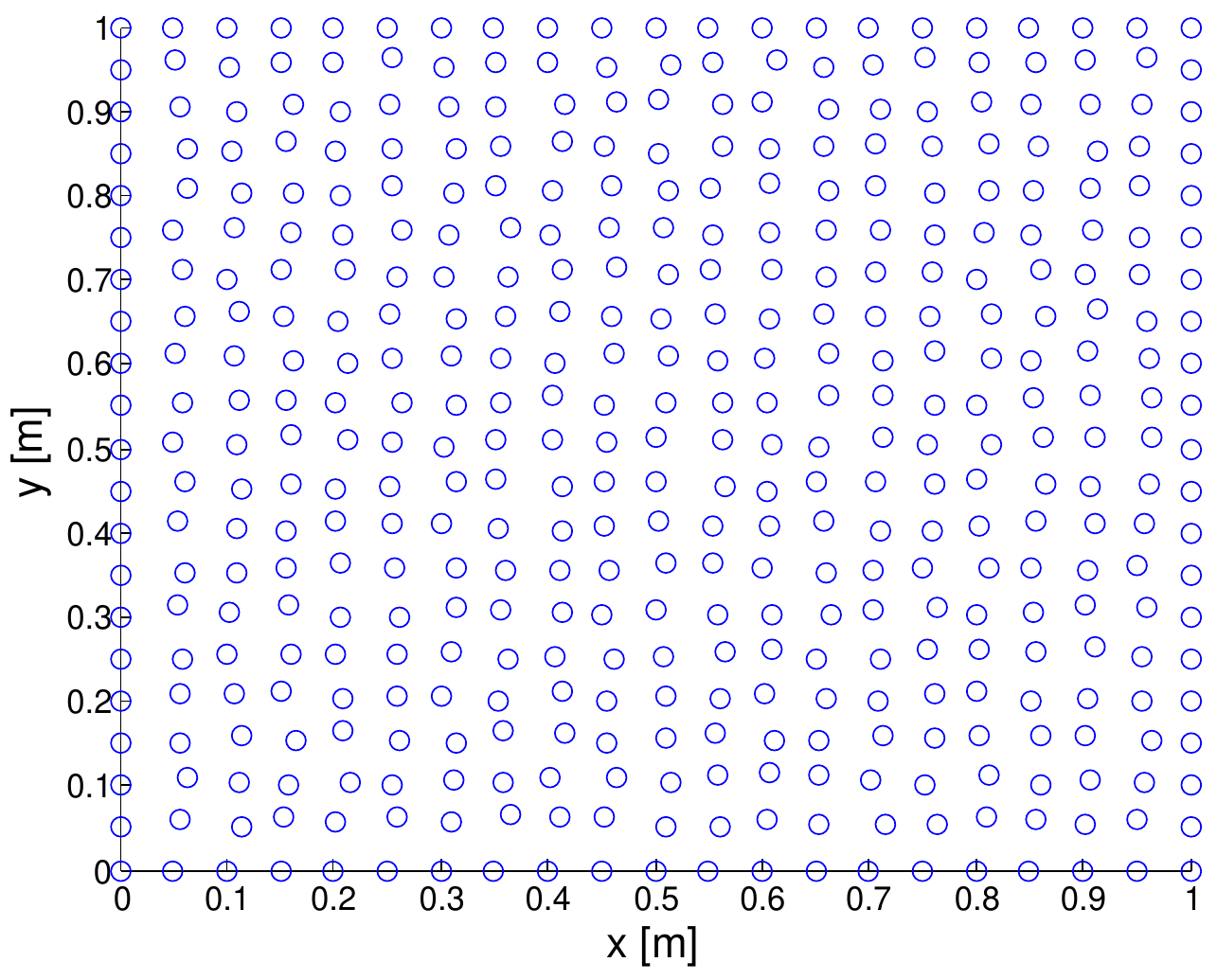}}\\
\caption{Irregular nodal distributions.}
\label{fig:P2mesh}
\end{figure}

\begin{figure}
 \centering
 \subfigure[regular node distribution]{
 \label{fig:P2hg}
 \includegraphics[width=.4\textwidth]{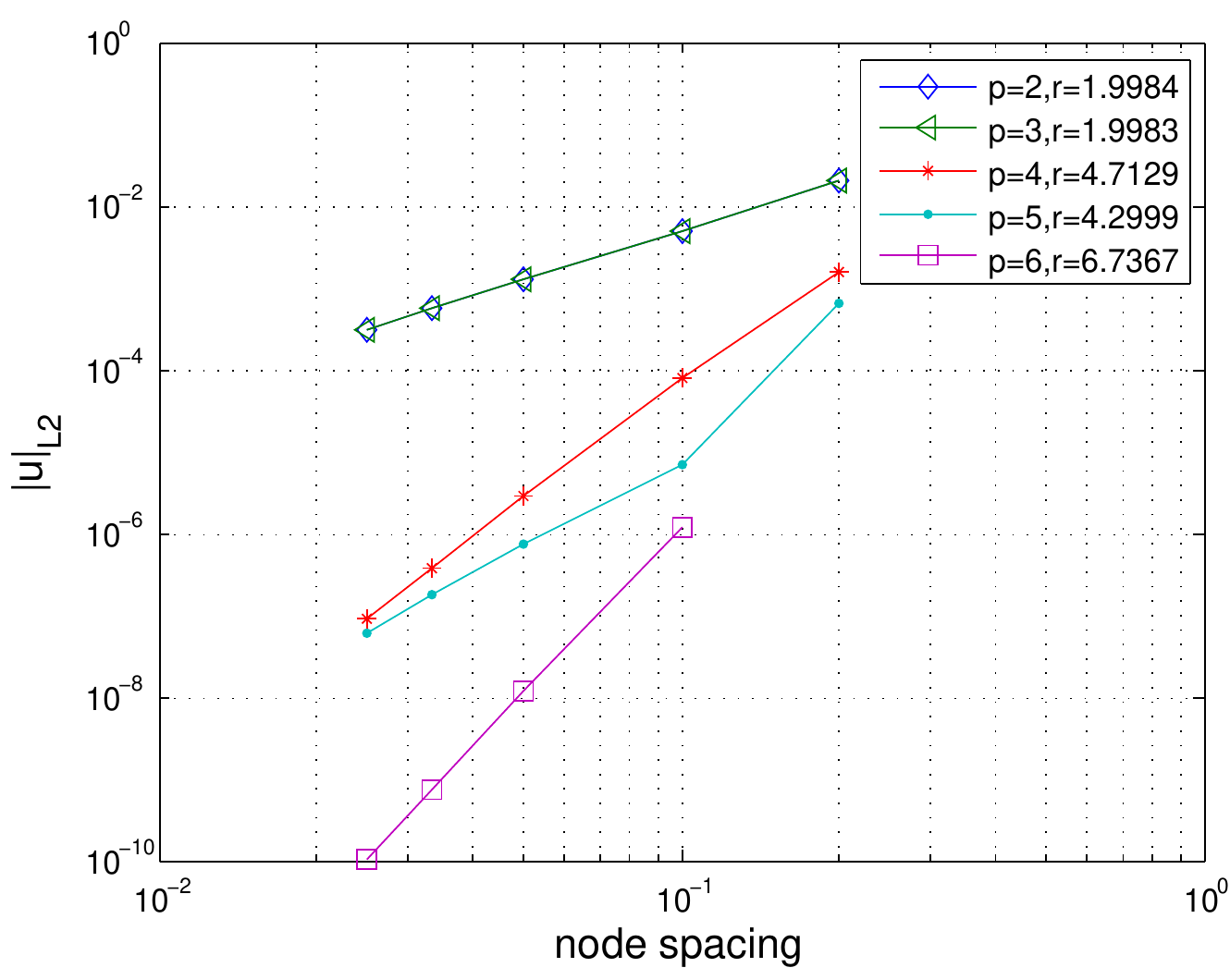}}
 \vspace{.1in}
 \subfigure[irregular node distribution]{
 \label{fig:P2noHg}
 \includegraphics[width=.4\textwidth]{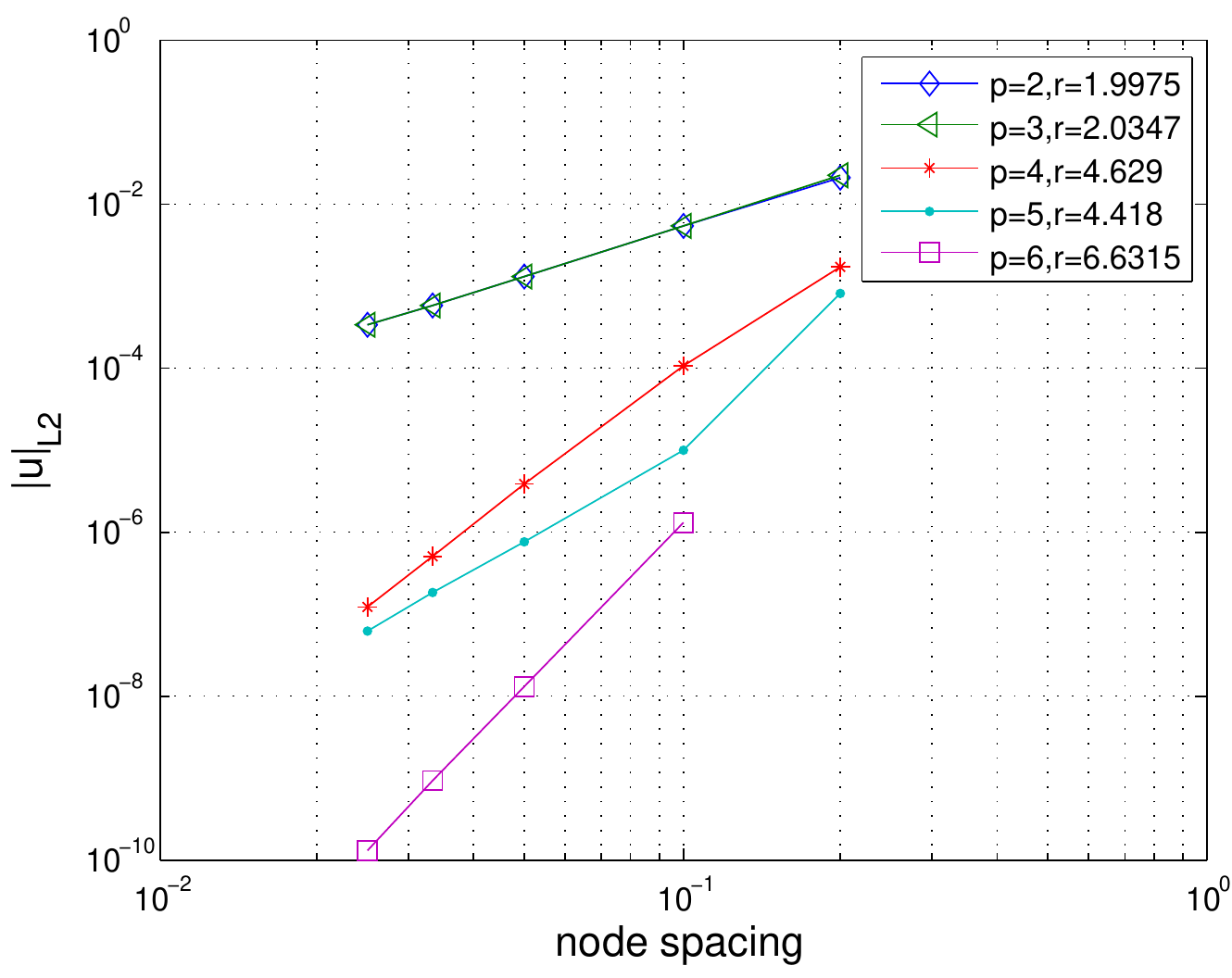}}\\
 \vspace{.3in}
\caption{The $L_2$ norm of different polynomial orders and node spacings for regular/irregular nodal distributions.}
\label{fig:P2contour}
\end{figure}

The number of neighbors for each point is selected as the number of terms in the nonlocal operator. We test the convergence of the L2 error for the $u$ field under uniform discretizations and non-uniform discretization in Fig.\ref{fig:P2mesh}. The convergent plot is given in Fig.\ref{fig:P2contour}.

\section{Numerical examples by weak form}\label{sec:nwf}
The fourth example aims at solving the Poisson equation in higher dimensional space by both the ``equivalent'' integral form and operator energy functional. The fifth example is about the biharmornic equation. The sixth example solves the Von-Karmon plate with simply support. 
\subsection{Poisson equation in higher dimensional space}
In this section, we solve the Poisson equation in $n$ dimensional space by nonlocal operator method. The $n$ dimensional Poisson equation is 
\begin{align}
\nabla^2 u=f(\bx),\quad \bx\in [0,1]^n \label{eq:npse}
\end{align}
with analytic solution
\begin{align}
u(\bx)=\exp{\Big(\sum_{i=1}^n (-1)^{i-1} x_i\Big)} \Pi_{i=1}^n x_i(1-x_i) 
\end{align}
under the boundary conditions
\begin{align}
u(x_1,...,x_i=0,...,x_n)=u(x_1,...,x_i=1,...,x_n)=0,\, 1\leq i\leq n
\end{align}
where $\bx=(x_1,...,x_n)$, $f(\bx)=\nabla^2 \big(\exp{(\sum_{i=1}^n (-1)^{i-1} x_i)} \Pi_{i=1}^n x_i(1-x_i)\big)$. 

The equivalent integral functional for Eq.\ref{eq:npse} is
\begin{align}
\fF=\int_{\Omega} \big(\frac{1}{2}\nabla u \cdot \nabla u-f(\bx) u\big) \ud V
\end{align}
The tangent stiffness is constructed from the operator matrix $\mathbf B$ for $\nabla u$, e.g.
\begin{align}
\mathbf K_g=\sum_{\Delta V_i \in \Omega} \Delta V_i \mathbf B_i^T \mathbf B_i
\end{align}
The Dirichlet boundary condition are applied by penalty method.

The Poisson equation with dimensional number $n=(2,3,4,5)$ under different discretization and order of   nonlocal operator are tested, where the statistical results are shown in Tables. (\ref{tab:modelABC2d},\ref{tab:modelABC3d},\ref{tab:modelABC4d},\ref{tab:modelABC5d}).
\begin{table}[htp]
\begin{center}
\begin{tabular}{cccccc}
\hline
Nnode & $\Delta x$ & $L_2$ norm & $\frac{u_{max}}{u_{max}^{exact}}-1$ & $p$-order & $p^{hg}$\\ \hline
1681& 0.025& 0.0485& -0.0281& 1& 1\\
1681& 0.025& 0.0262& 0.01& 2& 1\\
1681& 0.025& 0.0139& -0.00256& 3& 1\\
1681& 0.025& 0.0175& -0.00308& 4& 1\\
6561& 0.0125& 0.0379& 0.033& 1& 0\\
6561& 0.0125& 0.0179& 0.0714& 1& 1\\
6561& 0.0125& 0.011& 0.00505& 2& 1\\
25921& 0.00625& 0.0202& 0.0221& 1& 1\\
25921& 0.00625& 0.00501& 0.00266& 2& 1\\
25921& 0.00625& 0.00191& -0.000417& 3& 1\\
40401& 0.005& 0.00777& -0.00263& 1& 1\\
160801& 0.0025& 0.00291& 0.0007& 1& 1\\\hline
\end{tabular}
\caption{Statistical results for 2 dimensional Poisson equation}\label{tab:modelABC2d}
\end{center}
\end{table}

Table.\ref{tab:modelABC2d} gives the statistical results for 2D Poisson equation under different discretizations. When the order of nonlocal operator increased from 1 to 3, the L2 norm and error for $u_{max}$ decrease gradually as shown in several cases. However, for 4-order nonlocal operator, the result is not better than 3-order scheme. The 3-order scheme with 25921 nodes can achieve better result than 1-order scheme with 160801 nodes. The comparison between 5,6 rows shows that the operator energy functional has positive effect in improving the accuracy. In contrast with the scheme by strong form, the convergence property of weak form is slightly affected by the operator energy functional.

\begin{table}[htp]
\begin{center}
\begin{tabular}{cccccc}
\hline
Nnode & $\Delta x$ & $L_2$ norm & $\frac{u_{max}}{u_{max}^{exact}}-1$ & $p$-order & $p^{hg}$\\ \hline
10648& 0.04763& 0.0907& -0.0406& 1& 1\\
29791& 0.03333& 0.0604& -0.0248& 1& 1\\
68921& 0.025& 0.0485& -0.02& 1& 1\\\hline
\end{tabular}
\caption{Statistical results for 3 dimensional Poisson equation}\label{tab:modelABC3d}
\end{center}
\end{table}
For the 3D Poisson equation, we tested three cases with discretization ranged from 22, 31, 41 nodes in each direction. The statistical results are given in Table.\ref{tab:modelABC4d}. The L2 norm and error for $u_{max}$ decrease with the point grid space. When 41 nodes used for each direction, the L2 norm is approximately 5\%.
\begin{table}[htp]
\begin{center}
\begin{tabular}{cccccc}
\hline
Nnode & $\Delta x$ & $L_2$ norm & $\frac{u_{max}}{u_{max}^{exact}}-1$ & $p$-order & $p^{hg}$\\ \hline
14641& 0.1& 0.169& -0.0514& 1& 1\\
65536& 0.0667& 0.118& -0.0171& 1& 1\\
160000& 0.0526& 0.0983& -0.0203& 1& 1\\
810000& 0.0345& 0.0579& 0.00304& 1& 1\\
2560000&0.0256&0.0454&0.00152&1&1\\\hline
\end{tabular}
\caption{Statistical results for 4 dimensional Poisson equation}\label{tab:modelABC4d}
\end{center}
\end{table}
For the 4 dimensional Poisson equation, we tested four cases with discretization ranged from 11,16,20,30,40 nodes for each direction. The statistical results are given in Table.\ref{tab:modelABC4d}. The L2 norm and error for $u_{max}$ decrease with the point grid space. When 40 nodes used for each direction, the L2 norm is approximately 5\%.

\begin{table}[htp]
\begin{center}
\begin{tabular}{cccccc}
\hline
Nnode & $\Delta x$ & $L_2$ norm & $\frac{u_{max}}{u_{max}^{exact}}-1$ & $p$-order & $p^{hg}$\\ \hline
7776& 0.2& 0.229& -0.114& 1& 1\\
100000& 0.111& 0.181& -0.0944& 1& 1\\
1048576& 0.0667& 0.13& -0.0485& 1& 1\\
4084101& 0.05& 0.0985& -0.0352& 1& 1\\\hline
\end{tabular}
\caption{Statistical results for 5 dimensional Poisson equation}\label{tab:modelABC5d}
\end{center}
\end{table}
For 5 dimensional Poisson equation, when 16 nodes are assigned in each direction, the number of nodes reaches 1,048,576. More nodes in each direction will lead to the dimension disaster. The statistical results for different discretization are given in Table \ref{tab:modelABC5d}. The L2 norm and error for maximal $u$ decrease with the node spacing. We tested maximal 21 nodes in each direction (the computational scale is restricted by the computational power of a desktop PC), the L2 norm is approximately 9.85\% and the error for $u_{max}$ with respect to the theoretical solution is less than 4\%.
\subsection{Square plate with simple support}
The plate equation reads
\begin{align}
w_{,04}+2 w_{,22}+w_{,40}=\frac{q_0}{D_0},\quad (x,y)\in (0,1)\times (-1/2,1/2)\label{eq:nbh}
\end{align}
where $D_0=\frac{E t^3}{12(1-\nu^2)}$, with Dirichlet boundary conditions
\begin{align}
w(x,-1/2)=w(x,1/2)=0,\,x\in [0,1]\notag\\
w(0,y)=w(1,y)=0,\,y\in [-1/2,1/2]\notag
\end{align}
The analytic solution for the simply support square plate subjected to uniform load is denoted by \cite{timoshenko1959theory}
\begin{align}
w=\frac{4 q_0 a^4}{\pi^5 D_0}\sum_{m=1,3,...}^{\infty} \frac{1}{m^5}\Big(1-\frac{\alpha_m \tanh \alpha_m+2}{2 \cosh \alpha_m}\cosh\frac{2 \alpha_m y}{a}+\frac{\alpha_m}{2\cosh\alpha_m}\frac{2y}{a}\sinh\frac{2\alpha_m y}{a}\Big)\sin\frac{m \pi x}{a}
\end{align}
where $\alpha_m=\frac{m \pi}{2}$. 

The ``equivalent'' integral form for Eq.\ref{eq:nbh} is
\begin{align}
\fF_{plate}=\frac{1}{2}\tpartial w^T D_{plate} \tpartial w
\end{align}
where
\begin{align}
\tpartial w=(w_{yy},w_{xx},w_{xy})^T
\end{align}
\begin{align}
\mathbf D_{plate}=\frac{E t^3}{12(1-\nu^2)}\left[
\begin{array}{ccc}
 1 & \nu & 0 \\
 \nu & 1 & 0 \\
 0 & 0 & 2-2 \nu \\
\end{array}
\right]
\end{align}

The parameters for the plate include length $a=1$, thickness $t=0.01$ m and uniform pressure $q_0$=-100 N, Poisson ratio $\nu=0.3$, elastic modulus $E=30$ GPa and $D_0=\frac{E t^3}{12(1-\nu^2)}$

With the aid of nonlocal operator $\tpartial w$ and its operator matrix $\mathbf B$, the first and second variation of the energy functional are,
\begin{align}
\delta \fF_{plate}&=\sum_{\Delta V_i\in \Omega} \Delta V_i\big(\tpartial w_{i}^T \mathbf D_{plate} \tpartial \delta w_i-q_0 \delta w_i\big)=\sum_{\Delta V_i \in \Omega} \Delta V_i \tpartial w_{i}^T D_{plate} \mathbf B_i \delta \mathbf w_i-q_0 \delta w_i\notag\\
\delta^2 \fF_{plate}&=\sum_{\Delta V_i\in \Omega} \Delta V_i\tpartial \delta w_{i}^T \mathbf D_{plate} \tpartial \delta w_i=\sum_{\Delta V_i\in \Omega} \Delta V_i\delta \mathbf w_i^T \mathbf B_i^T \mathbf D_{plate} \mathbf B_i \delta \mathbf w_i.\notag
\end{align}
where $\delta \mathbf w_i$ is the vector for all unknowns in support $\cS_i$.

\begin{figure}
 \centering
 \subfigure[]{
 \label{fig:bhmW}
 \includegraphics[width=.4\textwidth]{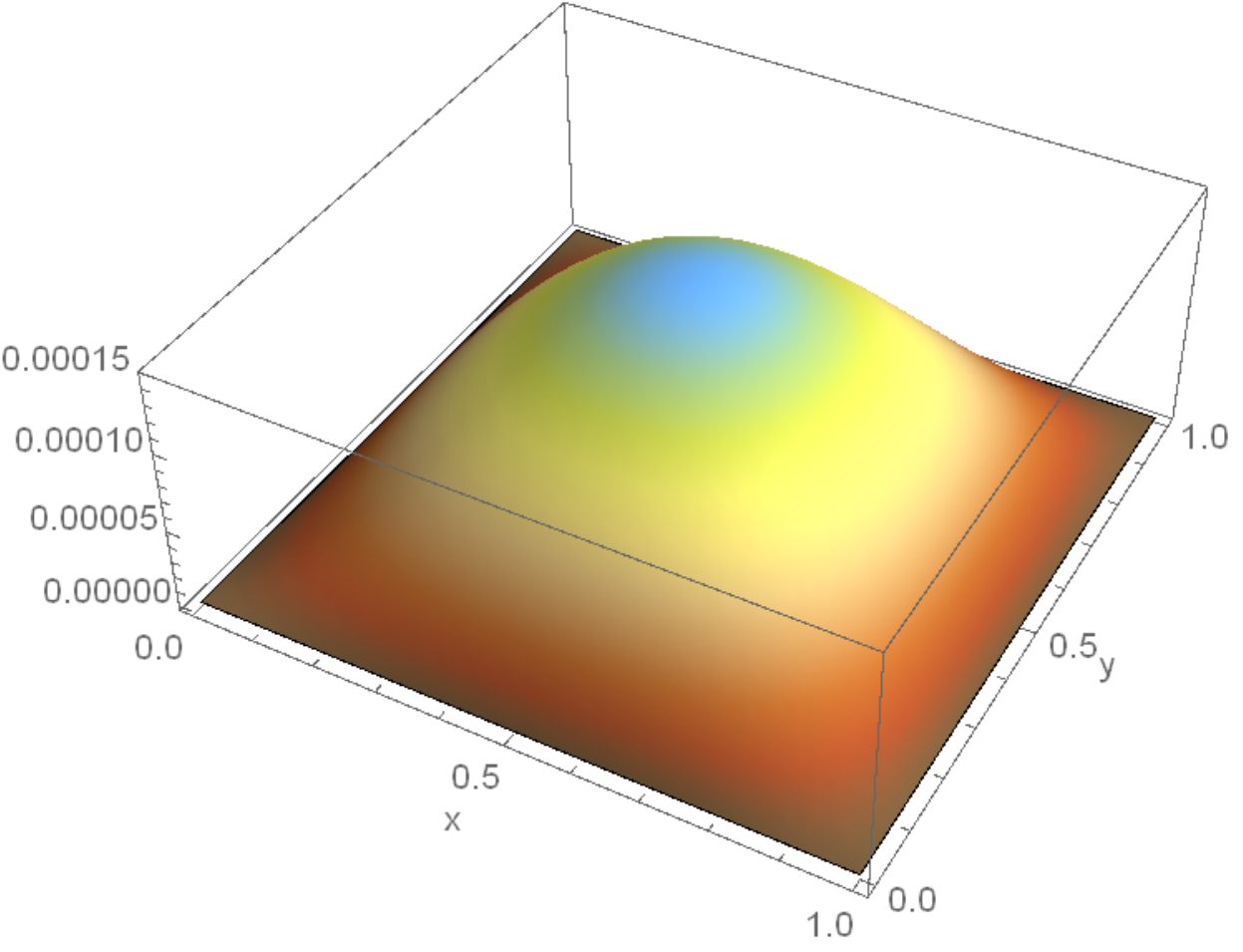}}
 \vspace{.1in}
 \subfigure[]{
 \label{fig:bhmC}
 \includegraphics[width=.4\textwidth]{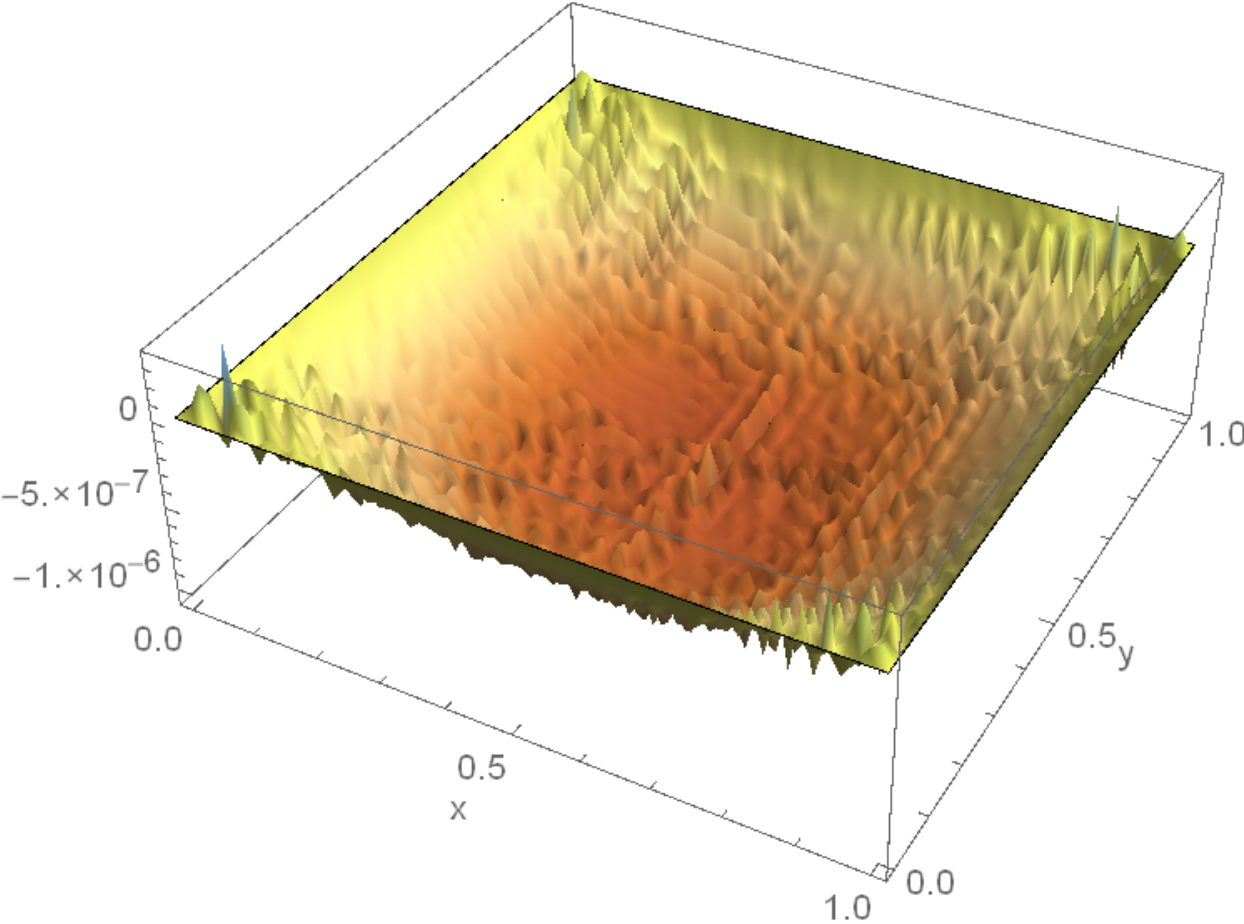}}\\
 \vspace{.3in}
\caption{(a) Deflection. (b) Error of deflection $w$ for discretization of 40$\times$40 with respect to exact solution.}
\label{fig:bhm}
\end{figure}

The plate is discretized uniformly and the number of neighbors for each point is selected as $n=5 p+length(\tpartial u)$, where $p$ is the order of nonlocal operator. The deflection curves for several discretizations are compared with the analytic solution in Fig.\ref{fig:bhmC}. The contour of the deflection field for discretization of $40\times 40$ is shown in Fig.\ref{fig:bhmC}. Compared with the original nonlocal operator method, the higher order NOM obtains the nonlocal operator in a simper way.

\begin{figure}
	\centering
		\includegraphics[width=6cm]{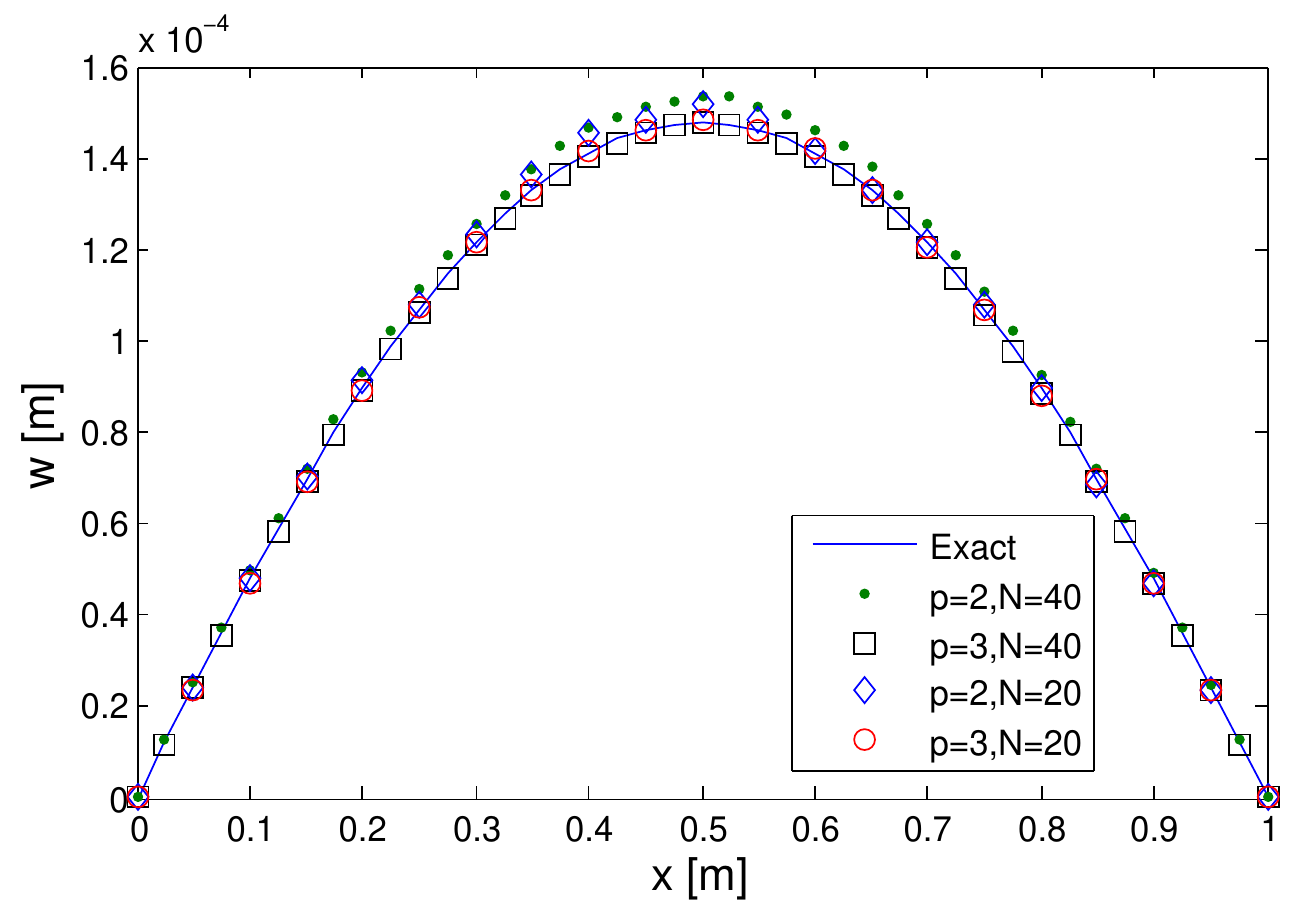}
\caption{Deflection of section $y=0$ under different discretizations, where $p$ denotes order of nonlocal operator, $N$ is the number of nodes in one direction.}\label{fig:bhmW0}
\end{figure}

\subsection{Von K\'arm\'a equations for a thin plate}
The Von K\'arm\'a equations \cite{von1910festigkeitsprobleme} are a set of nonlinear partial differential equations describing the large deflections of thin flat plates. The equations are based on Kirchhoff hypothesis : the surface normals to the plane of the plate remain perpendicular to the plate after deformation and  in-plane (membrane) displacements are small and the change in thickness of the plate is negligible. These assumptions imply that the displacement field $v$ in the plate can be expressed as \cite{ciarlet1980justification},
\begin{align}
 v_{1}(x_{1},x_{2},x_{3})&=u_{1}(x_{1},x_{2})-x_{3}\,{\frac {\partial w}{\partial x_{1}}},\notag\\
v_{2}(x_{1},x_{2},x_{3})&=u_{2}(x_{1},x_{2})-x_{3}\,{\frac {\partial w}{\partial x_{2}}},\notag\\
v_{3}(x_{1},x_{2},x_{3})&=w(x_{1},x_{2})
\end{align}
For a plate of a thickness $h$ defined on the mid-surface $(x_1,x_2)$, Von K\'arm\'a energy is given by \cite{landau1986theory,patricio1997numerical} 
\begin{align}
\iint_\Omega \{ \frac{E h^3}{24(1-\nu^2)} \{(\Delta w)^2-2(1-\nu) [w,w]\}+\frac{h}{2} \varepsilon_{ij} \sigma_{ij}-q w\}d x_1 d x_2.\label{eq:FVonKarman}
\end{align}
where Laplace operator $\Delta w=\frac{\partial^2 w}{\partial x_1^2}+\frac{\partial^2 w}{\partial x_2^2}$ and
\begin{align}
[w,w]=\frac{\partial^2 w}{\partial x_1^2} \frac{\partial^2 w}{\partial x_2^2}-(\frac{\partial^2 w}{\partial x_1 \partial x_2})^2
\end{align}
$\varepsilon_{ij}$ is the strain tensor with nonlinear terms in the deformations $u_1=u_1(x_1,x_2), u_2=u_2(x_1,x_2),w=w(x_1,x_2)$:
\begin{align}
\varepsilon_{ij}=\frac{1}{2}(\frac{\partial u_i}{\partial x_j}+ \frac{\partial u_j}{\partial x_i})+\frac{1}{2} \frac{\partial w}{\partial x_i}\frac{\partial w}{\partial x_j}
\end{align}
where $(u_1,u_2)$ is the lateral displacement field due to membrane effect, $w$ is the deflection, $\sigma_{ij}$ is the stress tensor, linearly proportional to $\varepsilon_{ij}$, $E=30\times 10^6$ Pa and $\nu=0.3$ are the Young modulus and the Poisson ratio, respectively, $q=1000 $Pa is the external normal force per unit area of the plate. The dimensions of the plate are $1.0 \times 1.0 \times 0.01 \mbox{ m}^3$. The energy functional in Eq.\ref{eq:FVonKarman} leads to the governing equations
\begin{align}
\frac{Eh^3}{12(1-\nu^2)}\nabla^4 w-h\frac{\partial}{\partial x_j}\left(\sigma_{ij}\frac{\partial w}{\partial x_i}\right)=q, \quad \frac{\partial\sigma_{ij}}{\partial x_j}=0
\end{align}
The Cauchy stress tensor in mid-plane can be written as
\begin{align}
\sigma&=\frac{E}{1-\nu^2} \big(\nu \mbox{tr} \varepsilon \mathbf I_{2\times 2}+(1-\nu) \varepsilon\big)\\
\varepsilon&=\begin{bmatrix}\varepsilon_{11} & \varepsilon_{12}\\
\varepsilon_{21}& \varepsilon_{22}\end{bmatrix}
\end{align}

In this paper, we write the moment and curvature by tensor form. The conventional vectorial form can be recovered with ease. The moment tenor and curvature tensor are
\begin{align}
\mathbf M=\begin{bmatrix}M_{11}& M_{12}\\ M_{21}& M_{22}\end{bmatrix}=D_0 \big(\nu \mbox{tr} \bm \kappa\bm I_{2\times 2}+ (1-\nu)\bm \kappa\big)
\end{align}
\begin{align}
\bm \kappa=\nabla \nabla w=\begin{bmatrix}\frac{\partial^2 w}{\partial x_1^2}&\frac{\partial^2 w}{\partial x_1\partial x_2}\\
\frac{\partial^2 w}{\partial x_2\partial x_1}& \frac{\partial^2 w}{\partial x_2^2} \end{bmatrix}
\end{align}
where $D_0=\frac{E h^3}{12 (1-\nu^2)}$. The moment tensor is similar to the stress tensor in the plane stress conditions.

The rotation in direction $\bm n$ is
\begin{align}
\frac{\partial w}{\partial n}=\nabla w \cdot \bm n, \mbox{ where }\bm n=(n_1,n_2)
\end{align}
The curvature in direction $\bm n$ is
\begin{align}
\kappa_n=\bm n^T \bm \kappa \bm n
\end{align}
The momentum in direction $\bm n$ is
\begin{align}
M_n=\bm n^T \mathbf M \bm n
\end{align}
The nonlocal differential operators in Eq.\ref{eq:FVonKarman} can be written as
\begin{align}
\tpartial u=(u_{1,01},u_{1,10},u_{2,01},u_{2,10},w_{,01},w_{,02},w_{,10},w_{,11},w_{,20})
\end{align}
The gradient of energy functional on $\tpartial u$ is
\begin{align}
\frac{\partial \fF}{\partial {\tpartial}\bu}={D_0}\begin{bmatrix}\frac{(1-\nu)}{2}  (u_{1,01} u_{2,10}+w_{,01} w_{,10})\\
\frac{1}{2} \left(\nu (2 u_{2,01}+w_{,01}^2)+2u_{1,10}+w_{,10}^2\right)\\
\frac{1}{2} \left(\nu (2u_{1,10}+w_{,10}^2)+2 u_{2,01}+w_{,01}^2\right)\\
\frac{(1-\nu)}{2}  (u_{1,01} u_{2,10}+w_{,01} w_{,10})\\
\frac{1}{2} \left( w_{,01}(2 \nu u_{1,10} +2 u_{2,01}+w_{,01}^2)+
 w_{,10}((1-\nu)(u_{1,01} u_{2,10})+w_{,01} w_{,10})\right)\\
	\frac{1}{12} h^2 (\nu w_{,20}+w_{,02})\\
	\frac{1}{2} \left((u_{1,01} u_{2,10})(w_{,01}-\nu w_{,01})+w_{,10}(2 \nu u_{2,01}+2
u_{1,10}+w_{,01}^2+w_{,10}^2)\right)\\
	\frac{1}{6} h^2 (1-\nu) w_{,11}\\
	\frac{1}{12} h^2 (\nu w_{,02}+w_{,20})
	\end{bmatrix}
 \end{align}
The Hessian matrix of $\fF$ can be obtained with ease by computing $\frac{\partial^2 \fF}{\partial {\tpartial}\bu^2}$. The solution can be obtained when using the Newton-Raphson method in \ref{sec:nrnom}. For simplicity, we only consider the simple support boundary conditions.

The plate solved by NOM is discretized by $50 \times 50 $ nodes. The reference results are calculated by S4R plate/shell element in ABAQUS \cite{hibbett1998abaqus}. S4R element is a 4-node doubly curved thin or thick shell element with reduced integration, hourglass control, finite membrane strains. In ABAQUS, the flat thin plate with the same material parameters are discretized into $100\times 100 $ elements. 

\begin{figure}
	\centering
		\includegraphics[width=10cm]{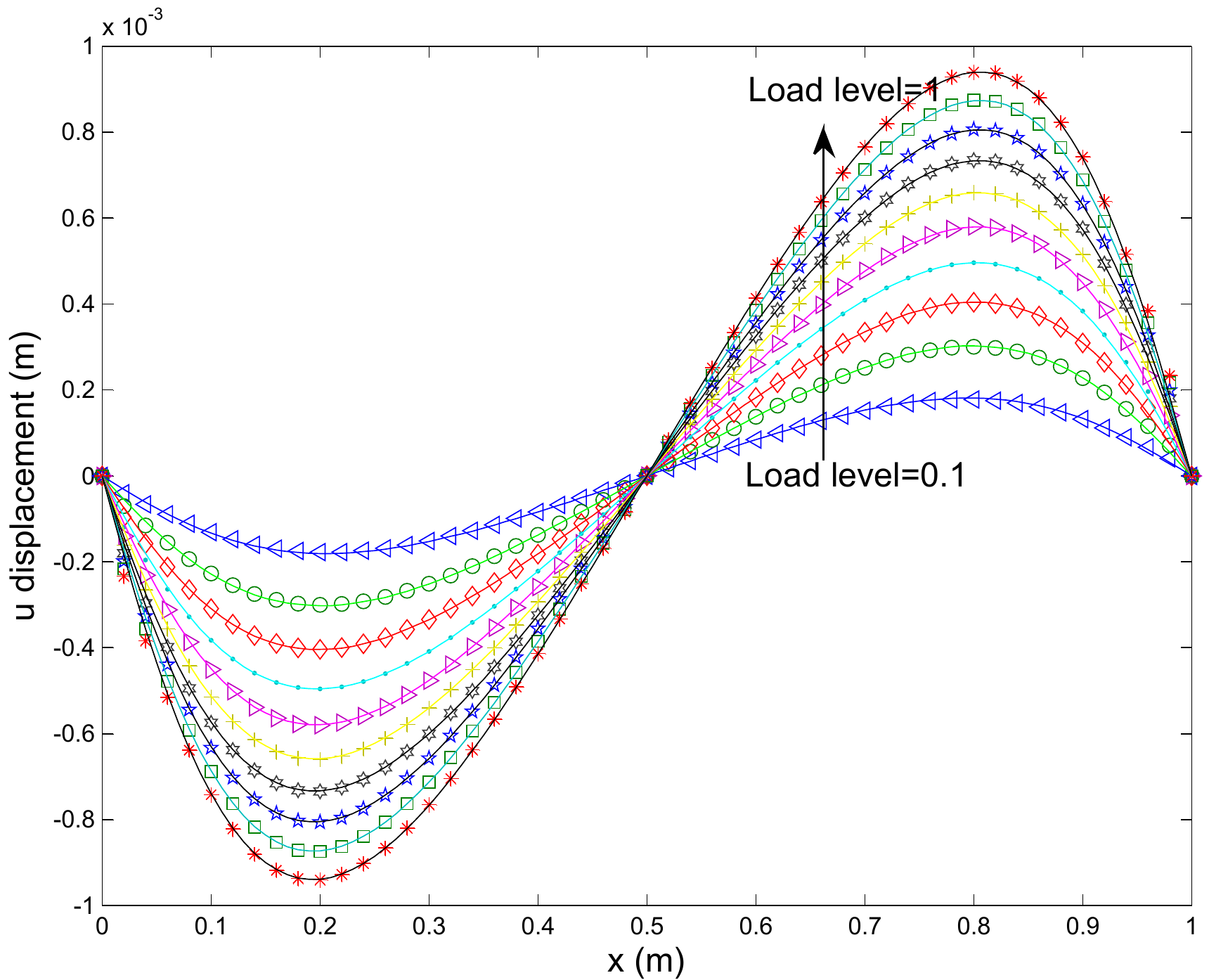}
	\caption{Displacement in membrane for nodes in $y=L/2$ under load level from 0.1 to 1, where the lines represent the numerical results by S4R element in ABAQUS while the star, diamond etc symbols are the results by NOM}
	\label{fig:VU2}
\end{figure}

\begin{figure}
	\centering
		\includegraphics[width=10cm]{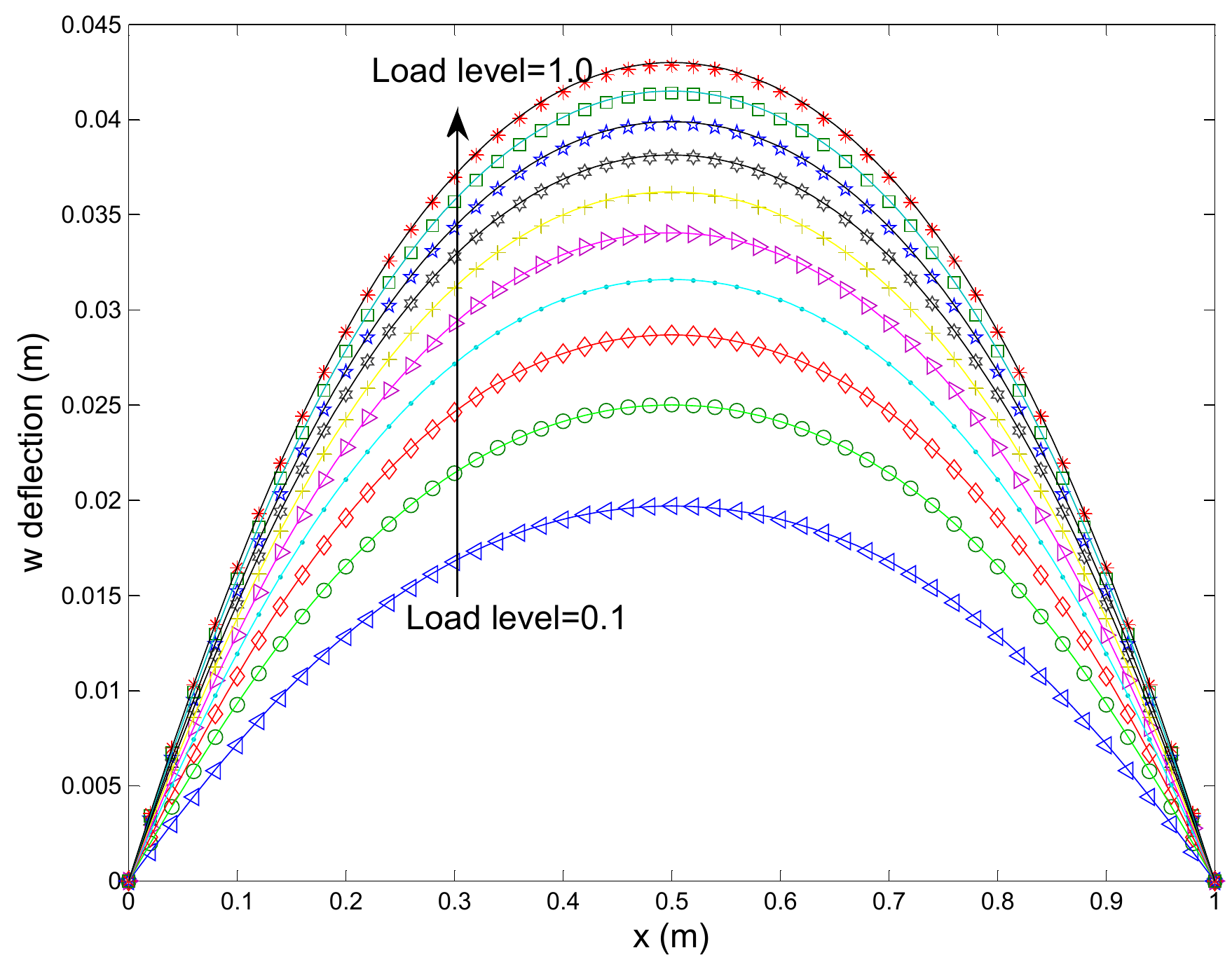}
	\caption{Deflection for nodes in $y=L/2$ under load level from 0.1 to 1, where the lines represent the numerical results by S4R element in ABAQUS while the star, diamond etc symbols are the results by NOM}
	\label{fig:VW2}
\end{figure}

Displacement in membrane and deflection out-of-plane for nodes on $y=0.5$ under different load levels are depicted in Figs.\ref{fig:VU2},\ref{fig:VW2}, respectively, where the lines represent the results by ABAQUS while the discrete symbols are the results by NOM. The displacement results agree well with that by ABAQUS.

\begin{figure}
	\centering
		\includegraphics[width=10cm]{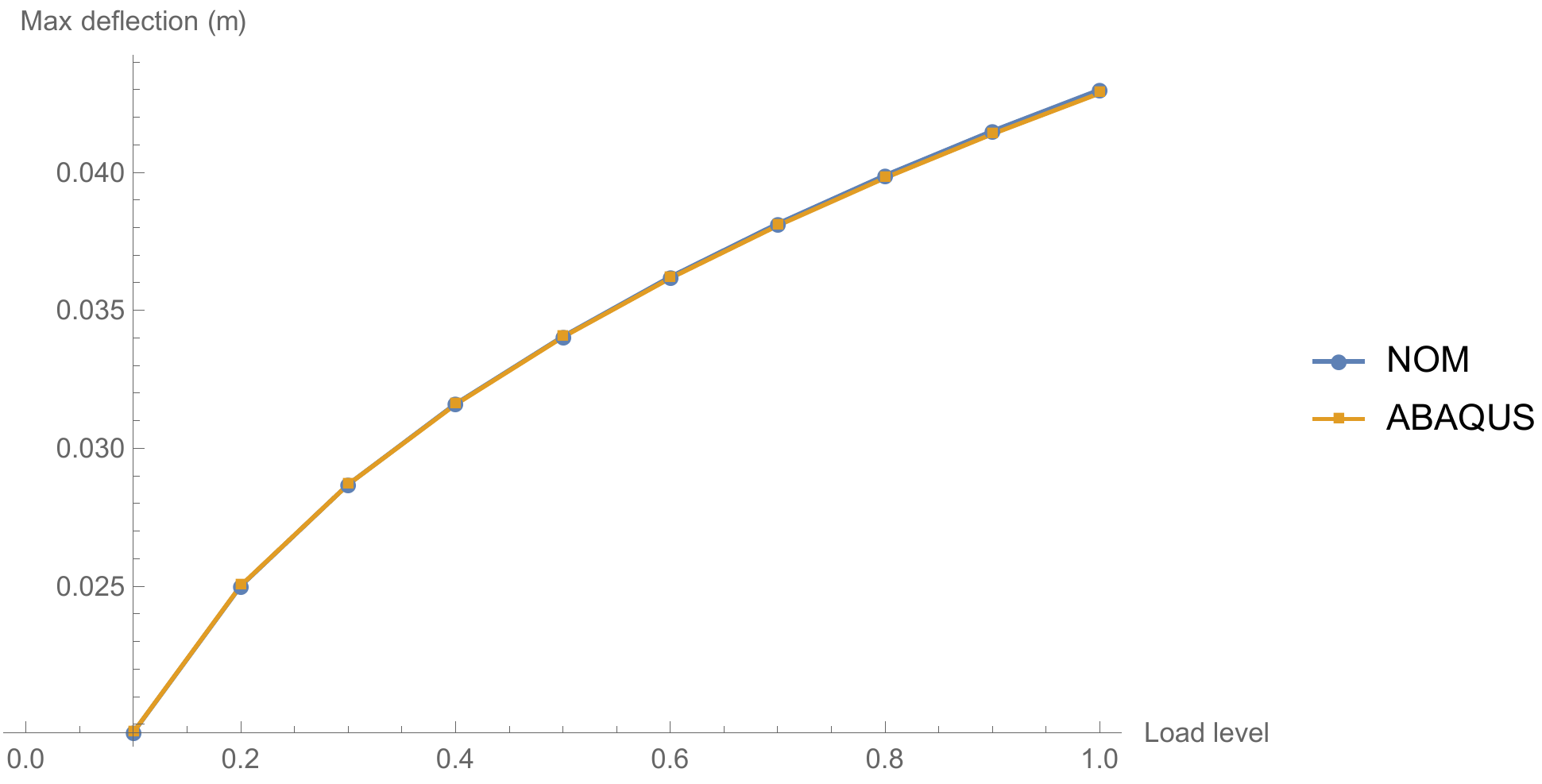}
	\caption{Maximal deflection for node in $(L/2,L/2)$ under load level from 0.1 to 1}
	\label{fig:VWmax}
\end{figure}

Maximal central deflection is plotted in Fig.\ref{fig:VWmax}, which shows the non-linearity increase with load level significantly. It can be seen that the result by NOM matches well with by finite element method.

\subsection{Nearly incompressible block}\label{sec:neoexample}
\begin{figure}
	\centering
		\includegraphics[width=5cm]{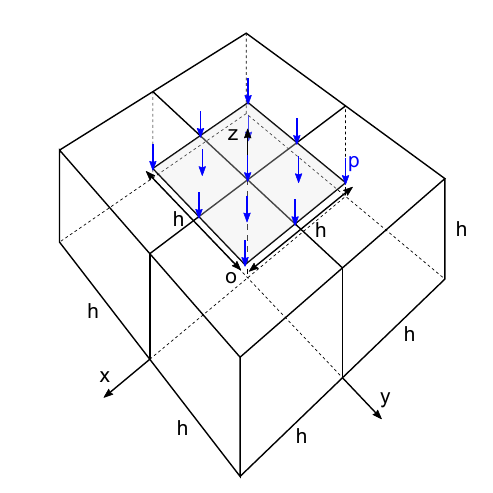}
	\caption{Setup of the block.}
	\label{fig:blockSetup}
\end{figure}

In this section, we model the nearly incompressible block of material consitituion in Eq.\ref{eq:nhm} by nonlocal operator method with Newton-Raphson iteration method. The nearly incompressible block of height $h=50$ mm, length $2 h$ and width $2h$ is loaded by an equally distributed pressure $p=3$ MPa at its top center of area $h \times h$ mm$^2$, as shown in Fig.\ref{fig:blockSetup}. For symmetry reason, only a quarter of the block is modeled. The bottom face is fixed in $z$-direction, while the nodes on plane $y=0$ are fixed in $y$-direction and the nodes on plane $x=0$ are fixed in $x$-direction, as similarly presented in reference \cite{reese2000new}. The material parameters are $\kappa=499.92568$ MPa, $\mu=1.61148$ MPa. 

\begin{figure}
 \centering
 \subfigure[FEM\cite{korelc2016automation} with $8^3$ mesh]{
 \label{fig:blockD1}
 \includegraphics[width=.22\textwidth]{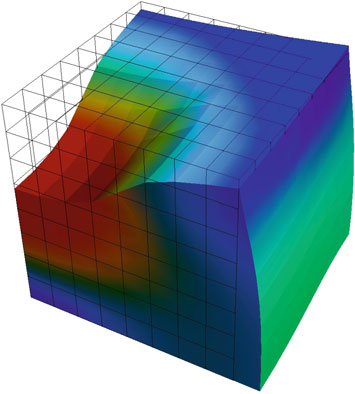}}
 \subfigure[NOM with $11^3$ nodes]{
 \label{fig:blockD2}
 \includegraphics[width=.3\textwidth]{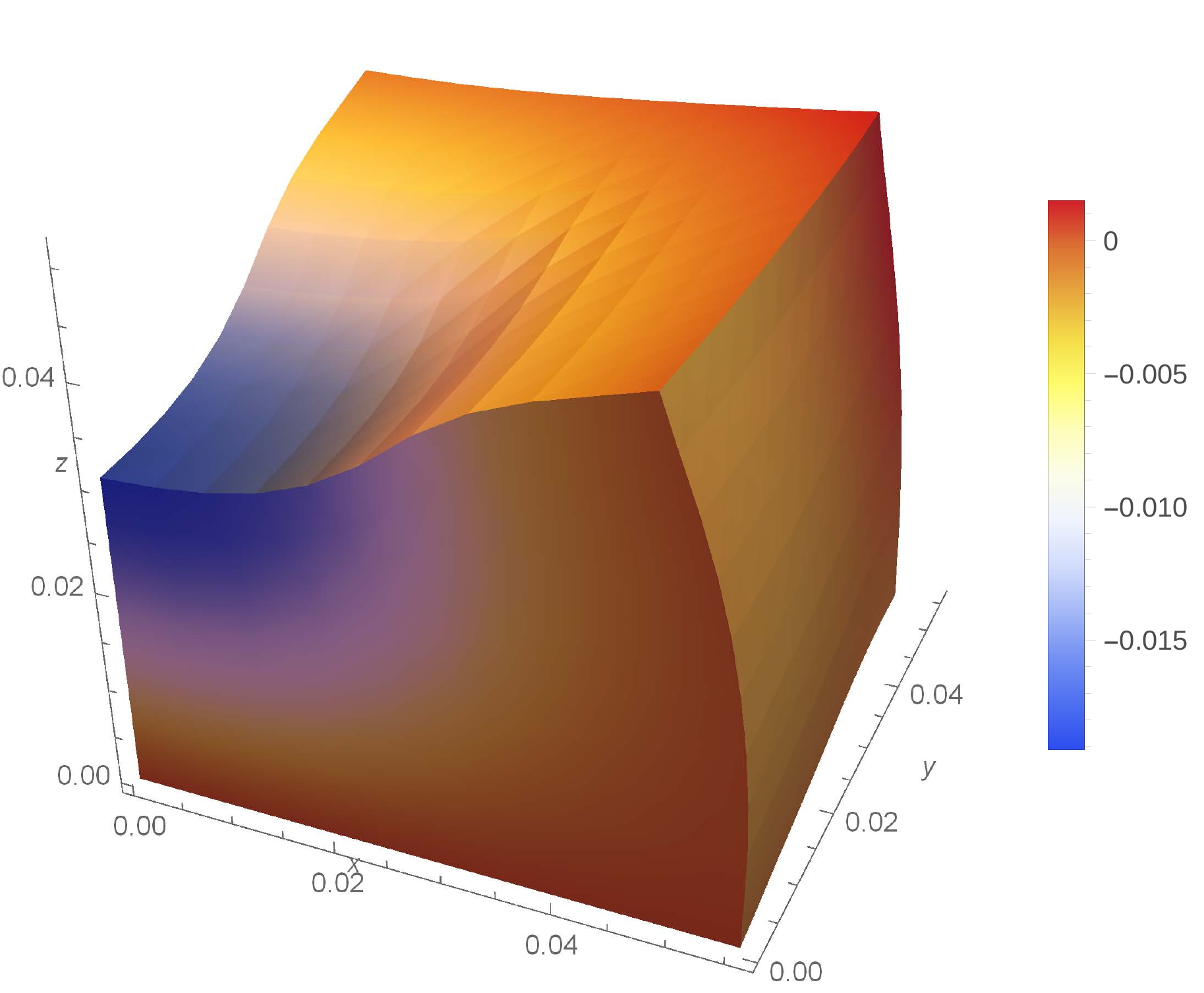}}
 \subfigure[NOM with $21^3$ nodes]{
 \label{fig:blockD3}
 \includegraphics[width=.33\textwidth]{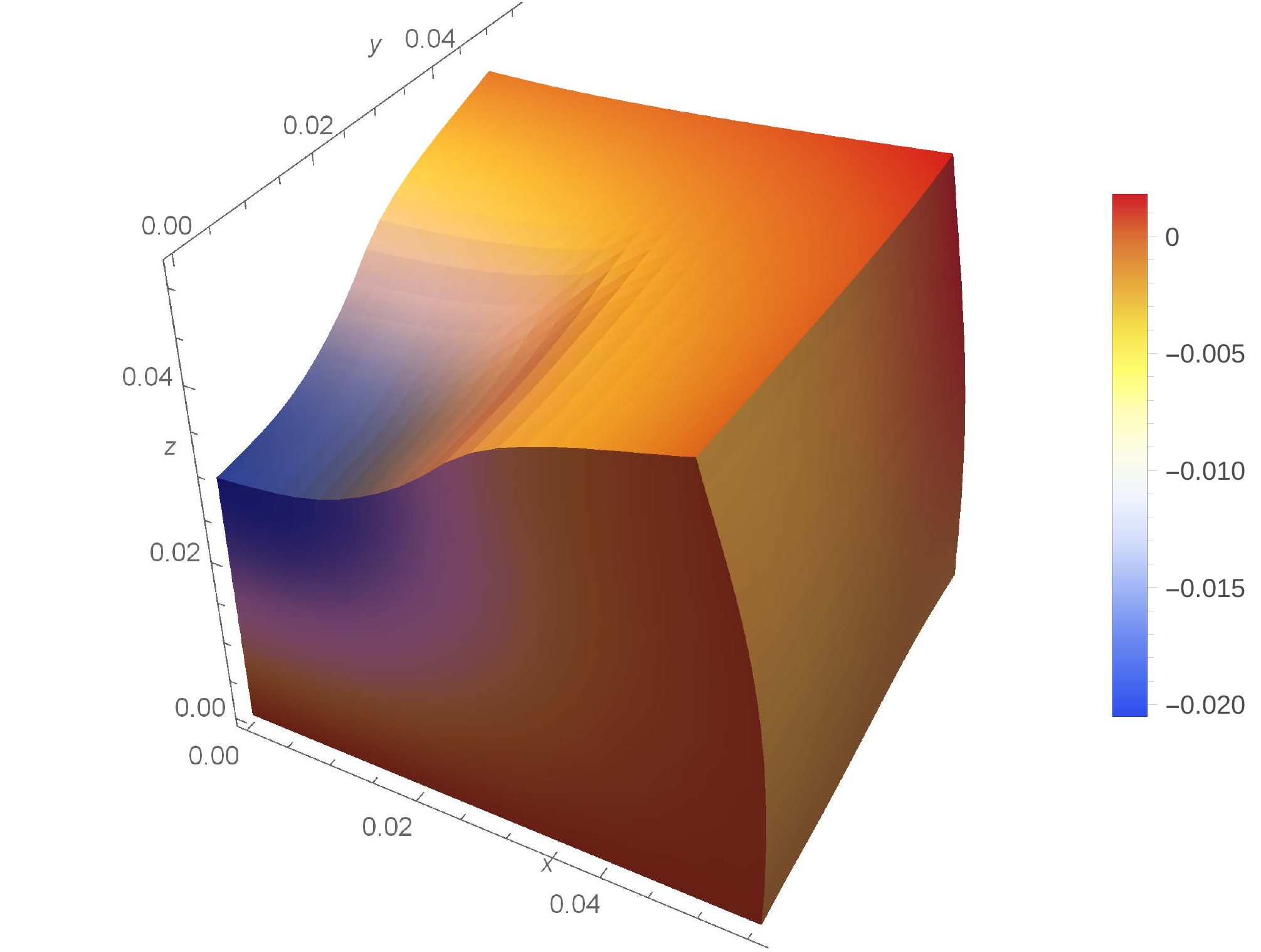}}
\caption{$z$-direction displacement in deformed configuration at final load level.}
\label{fig:blockD}
\end{figure}

The deformed block at final load level is depicted in Fig.\ref{fig:blockD}, where good agreement is obtained between the finite element method and the NOM. The maximal displacements in $z$-direction by linear hexahedral element (H1), quadratic hexahedral elements (H2) and nonlocal operator method are given in Table.\ref{tab:zmax}.

\begin{table}[h]
\begin{center}
\begin{tabular}{ | c|c | c |}
\hline
		H1 element &13.17 ($8^3$ mesh) &19.52 ($32^3$ mesh) \\ \hline
		H2 element &19.54 ($8^3$ mesh) &20.01 ($32^3$ mesh) \\ \hline
		NOM &19.14 ($11^3$ nodes) &20.43 ($21^3$ nodes)\\ \hline
\end{tabular}
\caption{Nearly incompressible block: displacement $w_{max}$ (mm)}\label{tab:zmax}
\end{center}
\end{table}

\section{Concluding Remarks}\label{sec:con}
\vspace{-2pt}
We have proposed a higher order nonlocal operator method for solving higher order PDEs based on the strong form or the equivalent integral formulated by weighed residual method or variational principles.

The relation of nonlocal operator and local operator is that the local operator is defined on a point, while the nonlocal operator method is defined on the support with finite characteristic length scale. When the support decreases to a point, the nonlocal operator degenerates to the local operator. Nonlocal operator is constructed from the Taylor series expansion and approximates the local derivative with orders up to $n$. In order to establish the nonlocal operator, only finite points in support is required. Nonlocal operator can be viewed as a generalization of the local operator. Most rules applied to the local operator can be adopted directly by the nonlocal operator method. 

In certain cases such as the regular grid, the nonlocal operator method is similar to the finite difference. One difference with finite difference method is that finite difference method requires a regular grid. When handling multiple fields, the finite difference method should adopt staggered grid for the reason of numerical stability, which complicates the numerical implementation. For nonlocal operator method, all the nodes have the same functions, in contrast with the finite difference method with staggered grid, where different nodes represent different fields. In terms of numerical stability, the nonlocal operator method monitors and enhances the robust of the derivative estimation by the operator energy functional, the quadratic functional of the Taylor series expansion. When adding the quadratic functional of the Taylor series expansion to the functional of physical problem, the numerical stability can be enhanced and traced.

Taylor series expansion of multiple variables based on the mutli-index notation is powerful in deriving various partial derivatives of different orders. Multi-index notation $\alpha_d^n$ in \ref{sec:mathcode} can obtain automatically all the partial derivatives with order up to $n$ in $d$ spatial dimensions. In addition, the characteristic length scale is introduced for high precision of the derivative estimation. With all the partial derivatives available, all linear PDEs up to $2n$ orders can be described with ease. By replacing the differential operator with the nonlocal one, nonlocal operator method converts the PDEs into algebraic equations directly. The nonlocal operator method can be viewed as a tool to study the higher order PDEs. 

\vspace{-6pt}
\section*{Acknowledgments}
The first author acknowledge the supports from the COMBAT Program (Computational Modeling and Design of Lithium-ion Batteries, Grant No.615132). The supports from National Basic Research Program of China (973 Program: 2011CB013800) and NSFC (51474157), the Ministry of Science and Technology of China (Grant No.SLDRCE14-B-28, SLDRCE14-B-31) are acknowledged.
\vspace{-3pt}
\appendix

\section{Taylor series expansion}\label{sec:Ty}
There are several formulations for the Taylor series expansion of function of multiple variables. 
The conventional Taylor series of a function at origin can be written as \cite{hormander1983analysis}
\begin{align}
u(x_{1},...,x_{d})&=\sum _{n_{1}=0}^{\infty }... \sum _{n_{d}=0}^{\infty }{\frac {x_{1}^{n_{1}}... x_{d}^{n_{d}}}{n_{1}!... n_{d}!}}\,\left({\frac {\partial ^{n_{1}+... +n_{d}}u}{\partial x_{1}^{n_{1}}... \partial x_{d}^{n_{d}}}}\right)(0,... ,0)\label{eq:ty1}\\&=u(0,... ,0)+\sum _{j=1}^{d}{\frac {\partial u(0,... ,0)}{\partial x_{j}}}x_{j}+{\frac {1}{2!}}\sum _{j=1}^{d}\sum _{k=1}^{d}{\frac {\partial ^{2}u(0,... ,0)}{\partial x_{j}\partial x_{k}}}x_{j}x_{k}+\notag\\&\qquad \qquad +{\frac {1}{3!}}\sum _{j=1}^{d}\sum _{k=1}^{d}\sum _{l=1}^{d}{\frac {\partial ^{3}u(0,... ,0)}{\partial x_{j}\partial x_{k}\partial x_{l}}}x_{j}x_{k}x_{l}+...\label{eq:ty2}
 \end{align}

By using the generalization of inner product, the Taylor series expansion is
\begin{align}
u_j=u_i+\nabla u_i \cdot \br+\frac{1}{2!} \nabla^2 u_i : \br^2+...+\frac{1}{n!} \nabla^n u_i \cdot^{(n)} \br^n+...\label{eq:ty3}
\end{align}
where $\br=\bx_j-\bx_i$, $\br^n=\br\otimes ...\otimes \br$, and $\cdot^{(n)}$ is the generalization of inner product, where two special cases are $\cdot^{(1)}=\cdot,\,\cdot^{(2)}=:$.

Or using the $d$-dimensional multi-index, the Taylor series expansion is
\begin{align}
u_j=\sum _{(\alpha_{1},...,\alpha_{d})\in \alpha}\frac {r_{1}^{\alpha_{1}}...r_{d}^{\alpha_{d}}}{\alpha_{1}!...\alpha_{d}!} u_{i,\alpha_1...\alpha_d}\label{eq:ty4}
\end{align}
where $\alpha = \{(\alpha_1,...,\alpha_d)|\alpha_i\in \mathbb N^0,1\leq i\leq d\}$. In this paper, Eq.\ref{eq:ty4} is adopt for Taylor series expansion.

\section{Mathematica code for multi-index}\label{sec:mathcode}
For the multi-indexes in 
\[
\alpha_d^n=\{(n_1,...,n_d)|1\leq\sum_{i=1}^d n_i\leq n,\, n_i\in \mathbb N^0, 1\leq i \leq d\},
\]
the Mathematica code with high efficiency is
\begin{verbatim}
MultiIndexList[d_,n_]:=Module[{a,b,c},a=Subsets[Range[d+n],{d}];
Do[c=a[[i]];b=c-1;b[[2;;]]-=c[[1;;-2]];a[[i]]=b,{i,Length[a]}];a[[2;;]]];
(*note: d=number of spatial dimensions, n=maximal order of derivative*)
\end{verbatim}

The number of elements in $\alpha_d^n$ can be determined by counting the combination of positive integer $k$ as a sum of $d$ non-negative integers up to non-commutativity. Imagine a line of $d+k-1$ positions, where each position can contain either a cat or a divider. If one has $k$ (nameless) cats and $d-1$ dividers, he can split the cats into $d$ groups by choosing positions for the dividers: $C_{k+d-1}^{d-1}=C_{k+d-1}^{k}$, where $C_n^d$ is binomial coefficient and can be written as $C_n^d=\binom{n}{d}=\frac{n!}{d!(n-d)!}$. The size of each group of cats corresponds to one of the non-negative integers in the sum.

Therefore, in $d$-dimension space, the number of  $k$-order derivatives by Eq.\ref{eq:ty3} is
\begin{align}
N(\nabla^k u_i)= C_{k+d-1}^{k}.
\end{align}
The number of all derivatives with maximal order $n$ in $d$ dimensional space is 
\[
N_d^n=\sum_{k=1}^n C_{k+d-1}^{k}=C_{n+d}^{n}-1
\]
In order to obtain the $n$-order derivatives, the number of neighbors in $\cS_i$ must be not less than $C_{n+d}^{n}-1$, so that the coefficient matrix for all derivatives is invertible. If more neighbors are in the support, the least square method can be used to find the approximation. The minimal number of neighbors in support is listed in Table.\ref{tab:NumSupport}.
\begin{table}[h]
\begin{center}
\begin{tabular}{ c|cccccc }
\hline
$N_d^n$  & n=1  & n=2 & n=3 & n=4 & n=5 & n=6\\ \hline
d=1 &1 &2 &3 &4 & 5& 6\\ 
d=2 & 2&5 &9 &14 & 20 &27\\ 
d=3 & 3& 9& 19& 34& 55 & 83\\ 
d=4 & 4& 14& 34& 69& 125 &209\\
d=5 & 5& 20& 55& 125& 251 &461\\
d=6 & 6& 27& 83& 209& 461 &923 \\\hline
\end{tabular}
\caption{Minimal number of neighbors in support. $d$=number of spatial dimensions, $n$=maximal order of derivatives}\label{tab:NumSupport}
\end{center}
\end{table}

\section{Newton-Raphson method for nonlinear functional}\label{sec:nrnom}
The core of NOM is the functional, which comprises with physical functional and the operator energy functional. The physical functional may contains the functional on the domain and other functional on the boundaries. In all,
\begin{align}
\fF(\bu)=\int_{\Omega} (\fF_1^{ph}(\bu)+\fF^{hg}(\bu)) \ud V+\int_{\partial \Omega} \fF_2^{ph}(\bu) d S
\end{align}
The first and second derivative on all unknowns lead to the residual and the tangent stiffness matrix, respectively
\begin{align}
\mathbf R&=\frac{\partial\fF}{\partial \bu}=\int_{\Omega} (\frac{\partial\fF_1^{ph}}{\partial \bu}+\frac{\partial\fF^{hg}}{\partial \bu}) \ud V+\int_{\partial \Omega} \frac{\partial\fF_2^{ph}}{\partial \bu} d S\\
\mathbf K&=\frac{\partial \mathbf R}{\partial \bu^T}=\frac{\partial^2\fF}{\partial \bu\partial \bu^T}=\int_{\Omega} (\frac{\partial^2\fF_1^{ph}}{\partial \bu\partial \bu^T}+\frac{\partial^2\fF^{hg}}{\partial \bu\partial \bu^T}) \ud V+\int_{\partial \Omega}\frac{\partial^2\fF_2^{ph}}{\partial \bu\partial \bu^T} d S
\end{align}

When any term in $\fF$ is nonlinear functional, the Newton-Raphson is required. The solution is updated by iteration in each step. In the $n$ step, the residual $\mathbf R( \mathbf u_n)=0$ is satisfied, $\mathbf R( \mathbf u_{n+1})$ in the next step can be approximated by Taylor series expansion
\begin{align}
\mathbf R( \mathbf u_{n+1})\approx \mathbf R( \mathbf u_{n})+\frac{\partial \mathbf R}{\partial \bu^T}|_{\bu=\mathbf u_{n}} \cdot (\mathbf u_{n+1}-\mathbf u_{n})
\end{align}
The solution in $n+1$ step can be obtained by the iterations 
\begin{align}
0=\mathbf R( \mathbf u^{k+1})\approx \mathbf R( \mathbf u^{k})+ \mathbf K( \mathbf u^{k}) \cdot \Delta \mathbf u^{k+1} \to \mathbf K(\mathbf u^k)\Delta \mathbf u^{k+1}=- \mathbf R(\mathbf u^k)
\end{align}
where $k$ denotes the iteration number in $n+1$ step, $\mathbf u^{0}=\mathbf u_n$, $\mathbf u^{k+1}=\mathbf u^k+\Delta \mathbf u^{k+1}$. When 
\[
\frac{\|\Delta \mathbf u^{k+1}\|}{\|\sum_{i=1}^{k+1} \Delta \mathbf u^i\|}\leq \mbox{Tol},
\]
the iteration converges.
\vspace{-6pt}
\section*{References}
\vspace{-3pt}
\bibliographystyle{ksfh_nat}
\bibliography{honom.bbl}

\end{document}